\documentclass[12pt]{amsart}
\usepackage{amsmath,amsthm,amssymb,euscript,times}

\usepackage[matrix,arrow,curve,cmtip]{xy}

\newdir{ >}{{}*!/-5pt/\dir{>}}             
\newdir{ |}{{}*!/-5pt/\dir{|}}             

\newtheorem{thm}{Theorem}[section]
\newtheorem{prop}[thm]{Proposition}
\newtheorem{cor}[thm]{Corollary}
\newtheorem{lemma}[thm]{Lemma}

\theoremstyle{definition}
\newtheorem{defn}[thm]{Definition}

\newtheorem{example}[thm]{Example}

\theoremstyle{remark}
\newtheorem{remark}[thm]{Remark}

\numberwithin{equation}{section}

\DeclareMathOperator{\Alg}{Alg}

\DeclareMathOperator{\card}{card}

\DeclareMathOperator{\cf}{cf}

\DeclareMathOperator{\colim}{colim\,}

\DeclareMathOperator{\Field}{Field}

\DeclareMathOperator{\fp}{fp}

\DeclareMathOperator{\graph}{Graph}

\DeclareMathOperator{\id}{id}

\DeclareMathOperator{\lingraph}{lin-Graph}

\DeclareMathOperator{\Mod}{Mod}

\DeclareMathOperator{\pres}{pres}

\DeclareMathOperator{\sh}{Sh}

\DeclareMathOperator{\Sub}{Sub}

\DeclareMathOperator{\Vect}{Vect}

\DeclareMathOperator{\wellgraph}{well-Graph}

\newcommand{\emb}{\mathsf{emb}}

\newcommand{\mono}{\mathsf{mono}}

\newcommand{\op}{{\text{op}}}

\newcommand{\rank}{\mathsf{rank}}

\newcommand{\Set}{\mathit{Set}}
\newcommand{\Ab}{\mathit{Ab}}

\newcommand{\str}{{\mathsf{Str}}}

\newcommand{\topos}{\EuScript}
\newcommand{\cate}{\mathcal}

\newcommand{\A}{\topos A}

\newcommand{\B}{\cate B}
\newcommand{\fB}{\topos B}

\newcommand{\C}{\cate C}

\newcommand{\D}{\cate D}

\newcommand{\F}{\topos F}

\newcommand{\gb}{\textbf{g}}

\newcommand{\K}{\cate K}

\newcommand{\LL}{\topos L}

\newcommand{\N}{\mathbb{N}}

\newcommand{\cQ}{\cate Q}
\newcommand{\QQ}{\mathbb{Q}}

\newcommand{\RR}{\mathbb{R}}

\newcommand{\s}{\textup{\bd{s}}}

\newcommand{\cS}{\cate S}

\newcommand{\X}{\cate X}

\newcommand{\ZZ}{\mathbb{Z}}

\newcommand{\W}{\cate W}

\newcommand{\nil}{\varnothing}				

\newcommand{\te}{\bd{\textup{1}}}

\newcommand{\lc}{\langle}					
\newcommand{\rc}{\rangle}					

\renewcommand{\leq}{\leqslant}				
\renewcommand{\geq}{\geqslant}				
\renewcommand{\preceq}{\preccurlyeq}		

\newcommand{\into}{\rightarrowtail}

\newcommand{\ra}{\rightarrow}
\newcommand{\lra}{\longrightarrow}
\newcommand{\la}{\leftarrow}

\newcommand{\inc}{\hookrightarrow}

\newcommand{\bd}{\textbf}

\newcommand{\llra}[1]{\stackrel{#1}{\lra}}	
\newcommand{\xra}[1]{\xrightarrow{#1}}		

\begin{document}

\title{Cellular objects and Shelah's singular compactness theorem}
\author[T. Beke and J. Rosick\'{y}]
{T. Beke$^*$ and J. Rosick\'{y}$^{**}$}
\thanks{$^*$ Supported by the NSA under Grant H98230-11-1-0172. $^{**}$ Supported by the Grant Agency of the Czech Republic under the grant P201/12/G028.}
\address{\newline T. Beke\newline
Department of Mathematics\newline
University of Massachusetts Lowell\newline
One University Avenue, Lowell, MA 01854, U.S.A.\newline
tibor\_beke@uml.edu
\newline J. Rosick\'{y}\newline
Department of Mathematics and Statistics\newline
Masaryk University, Faculty of Sciences\newline
Kotl\'{a}\v{r}sk\'{a} 2, 611 37 Brno, Czech Republic\newline
rosicky@math.muni.cz
}
\date{\today}
\begin{abstract}
The best-known version of Shelah's celebrated singular cardinal compactness theorem states that if the cardinality of an abelian group is singular, and all its subgroups of lesser cardinality are free, then the group itself is free.  The proof can be adapted to cover a number of analogous situations in the setting of non-abelian groups, modules, graph colorings, set transversals etc.  We give a single, structural statement of singular compactness that covers all examples in the literature that we are aware of.  A case of this formulation, singular compactness for cellular structures, is of special interest; it expresses a relative notion of freeness.  The proof of our functorial formulation is motivated by a paper of Hodges, based on a talk of Shelah.  The cellular formulation is new, and related to recent work in abstract homotopy theory.
\end{abstract}
\maketitle

\section*{Introduction}
The form of Shelah's celebrated Singular Compactness Theorem that is the oldest historically and prototypical in its formulation is the following: if $\mu$ is a singular cardinal and $A$ is an abelian group of size $\mu$ all of whose subgroups of cardinality less than $\mu$ are free, then $A$ itself is free.  In his breakthrough work~\cite{S}, however, Shelah already proved variants of singular compactness for algebras other than abelian groups, as well as for certain graph colorings and set transversals, to which Hodges~\cite{H} and Eklof and Mekler~\cite{EM} added numerous other examples.  These examples are unified, roughly speaking, by the fact that the proof works for them.  Both Shelah~\cite{S} and Hodges~\cite{H} (which is also based on Shelah's ideas)\ contain axiomatizations for the proofs to carry through.  The \textsl{intrinsic} meaning of singular compactness remains slightly mysterious.  The best formulation that seems to be available is the informal ``if $\mu$ is a singular cardinal and $S$ a structure all of whose substructures of cardinality less than $\mu$ are free, then $S$ itself is free.''  Cf.\ the title of Hodges~\cite{H} and the first sentence of Eklof~\cite{E}.

The goal of this article, essentially, is to make the above sentence precise: to formulate singular compactness by specifying how \textsl{structure} and \textsl{free} should be understood in this context.  The obvious guess, that structure should mean `first-order definable structure,' turns out to be too broad; and the guess that free should mean `free algebra on a set,' turns out to be too restrictive.

Our starting point is the straightforward translation of the prototypical case of singular compactness into the setting of category theory.  Consider the free abelian group functor $F:\Set\ra\Ab$.  An abelian group is free if and only if it is isomorphic to an object in the image of this functor.  Our desired form of the singular compactness theorem is: for suitable functors $F:\A\ra\fB$, if the size of an object $X$ of $\fB$ is singular, and every subobject of $X$ of lesser size is in the image of $F$, then $X$ itself is in the image of $F$.  (For the sake of brevity, we will write `is in the image of the functor' to mean `is isomorphic to an object in the image of the functor'.)

A category is \textsl{accessible} if it is equivalent to the category of models and homomorphisms of a set of sentences of the form $\phi\ra\psi$, where $\phi$ and $\psi$ belong to the positive-existential fragment of the logic $L_{\kappa,\kappa}$ for some $\kappa$.  This class of categories was identified by Makkai and Par\'{e}~\cite{MP} as having the right mix of properties to develop a categorical model theory of infinitary logics.  There exist language-independent characterizations of accessible categories, using only concepts of category theory.  One needs to assume that $\A$ is an accessible category with directed colimits, $\fB$ is a finitely accessible category (this is a sub-class of accessible categories with directed colimits, corresponding, roughly, to $\kappa=\omega$ in the above definition)\ and the functor $F$ preserves directed colimits.  These assumptions allow one to introduce a notion of \textsl{size} of objects that is determined purely by the ambient category; see Def.~\ref{obj-size}.  The assumptions on directed colimits are indispensable for creating transfinite chains of subobjects.  The final assumption is that \textsl{$F$-structures extend along morphisms}.  This is a simple diagrammatic condition; see Def.~\ref{extend}.  It is probably best thought of as expressing for categories and functors what the Steinitz exchange property expresses for matroids.  It captures that $F$ is a `free' functor and $\A$ is a category of `bases' whose morphisms are `extensions of bases'.

\textsl{Singular compactness theorem (functorial form):}  Let $\A$ be an accessible category with filtered colimits, $\fB$ a finitely accessible category and $F:\A\ra\fB$ a functor preserving filtered colimits.  Assume that $F$-structures extend along morphisms.  Let $X\in\fB$ be an object whose size $\mu$ is a singular cardinal.  If all subobjects of $X$ of size less than $\mu$ are in the image of $F$, then $X$ itself is in the image of $F$.

The formulation given in the paper, Theorem~\ref{main}, is slightly more general in that it only assumes that `enough' subobjects of $X$ lie in the image of $F$.  The actual criterion, using dense filters of subobjects, was inspired by the treatment of singular compactness in Eklof-Mekler~\cite{EM}.

Returning to the paradigmatic example of singular compactness, there is another way to think of free abelian groups or more generally, free algebras.  Define a class $F_\alpha$ of algebras, $\alpha$ ranging over the ordinals, and compatible morphisms $F_\beta\ra F_\alpha$ for $\beta\prec\alpha$, by transfinite induction.  Let $F_\nil$ be $F(\nil)$, the free algebra on an empty set, and let $F_\bullet$ be the free algebra on a singleton.  For successor $\alpha+1$, let $F_{\alpha+1}$ be the pushout
\[  \xymatrix{ F_\nil \ar[r]\ar[d] &  F_\alpha \ar[d] \\
               F_\bullet \ar[r] & F_{\alpha+1}  } \]
For limit $\alpha$, let $F_\alpha$ be $\colim_{\beta\prec\alpha} F_\beta$.  Then an algebra is free if and only if it is isomorphic to $F_\alpha$ for some $\alpha$.

The above diagram is really a coproduct, since $F_\nil$ is the initial object of $\Alg$.  Working in an arbitrary cocomplete category and allowing an arbitrary member of a fixed collection $I$ of morphisms to be pushed on at successor ordinals, the process being continuous at limit ordinals, one obtains the important class of \textsl{$I$-cellular} morphisms.  An object $X$ is then called $I$-cellular if the (unique)\ map from the initial object to $X$ is $I$-cellular.

These concepts were introduced originally in topology, in the study of cellular spaces, where the inclusion $\partial D^n\inc D^n$ of the boundary spheres of the unit ball in $\RR^n$ (for any $n\in\N$)\ played the role that the morphism $F_\nil\ra F_\bullet$ did in the above example.  Classical algebraic topology is concerned with cell complexes built up in finitely many steps, but subsequent work extended this study into the transfinite.  See \cite{MRV} or the monographs of Hirschhorn~\cite{HM} or Hovey~\cite{HO} for background on cellular morphisms in an axiomatic setting and their relation to homotopy theory, and Def.~\ref{cell-def} for detailed definitions.

\textsl{Singular compactness theorem (cellular form):} Let $\fB$ be a locally finitely presentable category and $I$ a set of morphisms.  Let $X\in\fB$ be an object whose size $\mu$ is a singular cardinal.  If all subobjects of $X$ of size less than $\mu$ are $I$-cellular, then so is $X$.

Perhaps this is the most elegant version of singular compactness.  It can be thought of as formulating a relative version of freeness; `free structures' correspond to `free cell complexes,' that is, coproducts of cells.  The cellular form, Theorem~\ref{cellular}, covers all but a couple of classical examples of singular compactness that we are aware of.  It can be understood, for example, as a non-additive version of Eklof's formulation \cite{E} where the fat small object argument from \cite{MRV} replaces the Hill lemma.  We discovered this form while analyzing Shelah's notion of \textsl{$\mu$-colorable graph}: a graph whose vertices can be well-ordered so that no vertex is connected to more than $\mu$ of its predecessors in the well-ordering.  This well-ordering is just a trace of a cellular presentation; see Example~\ref{wellgraph}.

Cellular singular compactness is a corollary of functorial singular compactness.  Just as an abelian group has many bases, a cellular object has many cellular presentations, but one can code cells, attaching and cellular maps by data that amount to `bases' and `extensions of bases,' allowing one to reduce Theorem.~\ref{cellular} to Theorem~\ref{main}.  The requisite combinatorial structure was first identified, in a different setting, by Lurie~\cite{L}.  The verification of the basis extension property, condition (iii)\ of Theorem~\ref{main}, is non-trivial and builds on recent work in \cite{MRV}.

Our proof of functorial singular compactness is motivated by a close reading of Hodges~\cite{H}, itself based on a talk of Shelah's.  As with several papers of Hodges, there are ideas in \cite{H}, just beneath the surface, that can be readily phrased in the language of category theory, but the eventual form of our main result bears no easily discernible relation to the main result of \cite{H}.  We discuss this interaction more fully at the end of section~\ref{bigpict}.

\section{The main theorem}

Let $\A$, $\fB$ be categories and $F:\A\ra\fB$ a functor.  The following property will play a key role:

\begin{defn}   \label{extend}
We say that \textsl{$F$-structures extend along morphisms} if, given any morphism $g:X\ra Y$ and object $U$ of $\A$, together with an isomorphism $i:F(U)\ra F(X)$ in $\fB$, there exists a morphism $f:U\ra V$ and isomorphism $j:F(V)\ra F(Y)$ such that
\[ \xymatrix{  F(U) \ar[r]^{F(f)}\ar[d]^i &  F(V)\ar[d]^j \\
 F(X) \ar[r]^{F(g)} & F(Y) } \]
commutes.
\end{defn}

There is an equivalent, and perhaps slightly more suggestive, way of thinking of this property.  Throughout this paper, ``image'' means ``image up to isomorphism''; we will say that an object $X\in\fB$ is in the image of $F$ if there exists $U\in\A$ such that $F(U)$ is isomorphic to $X$, and we will say that a morphism $g:X\ra Y\in\fB$ is in the image of $F$ if there exist $f:U\ra V$ and isomorphisms $i,j$ such that the diagram
\[ \xymatrix{  F(U) \ar[r]^{F(f)}\ar[d]^i &  F(V)\ar[d]^j \\
 X \ar[r]^g & Y } \]
commutes.  

Then $F$-structures extend along morphisms if and only if for every $g:X\ra Y$ in the image of $F$ and isomorphism $i:F(U)\ra X$, there exist $f:U\ra V$ and isomorphism $j:F(V)\ra Y$ making the above diagram commute.  Thinking of an isomorphism $F(U)\ra X$ as an `$F$-structure on $X$', this condition states that if a morphism $g$ of $\fB$ underlies a morphism of $F$-structures, then any `reparametrization' of the $F$-structure of the domain can be extended to a reparametrization of the $F$-structure of $g$.

\begin{example}
Let $V$ be a finite-dimensional vector space and $\A$ the set of linearly independent subsets of $V$.  Order $\A$ by inclusion, turning it into a poset (hence category).  Let $\fB$ be the set of linear subspaces of $V$, again ordered by inclusion.  Let the functor (i.e.\ monotone map) $F:\A\ra\fB$ send a set of vectors to its linear span.  That $F$-structures extend is implied by the Steinitz exchange lemma.
\end{example}

For any object $X$ of a category $\C$, let $\Sub(X)$ denote the full subcategory of $\C/X$ whose objects are monomorphisms.  As usual, we will often denote an element $U\ra X$ of $\Sub(X)$ as $U$; this will cause no confusion.  For $U,V\in\Sub(X)$, we will say that $U$ \textsl{contains} $V$ and write $V\into U$ if there is a morphism (necessarily mono and unique)\ $V\ra U$ in $\Sub(X)$.

Let $\fB$ be an accessible category with filtered colimits, $X$ an object of $\fB$ and $\kappa$ a regular cardinal.

\begin{defn}   \label{filter}
A \textsl{$\kappa$-filter of $\Sub(X)$} is a full subcategory $\F\subseteq\Sub(X)$ such that for every directed system $\D\ra\F$ with $\card(D)\leq\kappa$ taking values in objects $U\in\F$ with $U$ $\kappa$-presentable, the colimit of $\D$ exists and belongs to $\F$.

A $\kappa$-filter $\F$ is said to be \textsl{dense} in $\Sub(X)$ if for every $\kappa$-presentable $S\in\Sub(X)$ there exists a $\kappa$-presentable $U\in\F$ with $S\into U$.

Let now $F:\A\ra\fB$ be a functor.  We say that the image of $F$ \textsl{contains a $\kappa$-dense filter of} $\Sub(X)$ if there exists a $\kappa$-dense filter $\F\subseteq\Sub(X)$ such that every $U\in\F$ is in the image of $F$.
\end{defn}

In any accessible category $\C$, for an object $X$ one can define
\[   \rank(X) = \min\{\kappa\;|\;X\textup{\ is\ }\kappa\textup{-presentable}\}   \]
If $\C$ is $\lambda$-accessible, has all filtered colimits and $\rank(X)>\lambda$ then $\rank(X)$ is a successor cardinal.  (See \cite{BR}~4.2). Thus we can introduce

\begin{defn}   \label{obj-size}
Let $\C$ be an accessible category with filtered colimits and $X\in\C$.  If $\rank(X)>\lambda$ where $\lambda$ is the least cardinal such that $\C$ is $\lambda$-accessible then denote by $\|X\|$ the cardinal predecessor of $\rank(X)$.
\end{defn}

For example, for $X\in\Set$, $\|X\|=\card(X)$ if $X$ is infinite, and undefined if $X$ is finite.  More generally, if $\C$ is the category of models and elementary embeddings, or category of models and homomorphisms of suitable (multi-sorted)\ theories in the logic $L_{\kappa,\lambda}$, then $\|X\|$ is the cardinality of the (disjoint union of the)\ sets underlying the model, provided it is no smaller than the L\"{o}wenheim-Skolem number of the theory.  $\|X\|$ depends, however, only on the ambient category $\C$, and not on the way it is axiomatized.  See \cite{BR}, 2.3 for further discussion.  We will use this concept rather than presentability when the `exact size' of the object needs to be captured.

\begin{defn}  \label{func-rank}
Let $F:\A\ra\fB$ be a functor preserving filtered colimits between accessible categories with filtered colimits.  The \textsl{rank of accessibility} of $F$ is the least cardinal $\mu_F$ such that $\A$ and $\fB$ are $\mu_F$-accessible, and for all regular $\lambda\geq\mu_F$, $F$ takes $\lambda$-presentable objects to $\lambda$-presentable ones.
\end{defn}

It is a basic fact about accessible categories that $\mu_F$ exists; see Prop.~\ref{unif}.

For a finitely accessible category $\fB$, let $\fp\fB$ denote the full subcategory of $\fB$ formed by one representative from each isomorphism class of its finitely presentable objects.  For any small category $\C$, let $\card(\C)$ denote the greater of $\aleph_0$ and the cardinality of the set of morphisms of $\C$.  For finitely accessible $\fB$, one can think of $\card(\fp\fB)$ as the size of a certain, canonically definable, first-order theory whose category of models is equivalent to $\fB$; Prop.~\ref{main1} will make this clear.  Now we can state the main result of this paper.

\begin{thm}   \label{main}
Let $\A$ be an accessible category with filtered colimits, $\fB$ a finitely accessible category and $F:\A\ra\fB$ a functor preserving filtered colimits.  Let $X\in\fB$ be an object with $\max\{\mu_F,\card(\fp\fB)\}<\|X\|$.  Assume
\begin{itemize}
\item[(i)] $\|X\|$ is a singular cardinal
\item[(ii)] there exists $\phi_{F,X}<\|X\|$ such that for all regular $\kappa$ with $\phi_{F,X}<\kappa<\|X\|$, the image of $F$ contains a dense $\kappa$-filter of $\Sub(X)$
\item[(iii)] $F$-structures extend along morphisms.
\end{itemize}
Then $X$ is in the image of $F$.
\end{thm}

The next section collects some preliminary facts about accessible categories.  This is followed by a section on two-player games for constructing countable chains of subobjects, a beautiful idea of Shelah's that is at the heart of all `simple' proofs of his singular cardinal compactness theorem.

The proof of the main theorem is in two stages.  In the first one, the target category $\fB$ is taken to be a presheaf category, i.e.\ the category of functors from some small category to $\Set$.  This will allow us to work with elements, making the proof much more transparent.  The second stage will extend the conclusion to any target category that possesses a functor into a presheaf category satisfying certain properties.  All finitely accessible categories do, as do all categories arising from Shelah's Abstract Elementary Classes.  See the last part of section 4 for a discussion about the generality of the main theorem.

Section 5 introduces cellular constructions and well-founded diagrams parametrizing them, and proves the cellular version of singular compactness as a corollary of the main theorem.  We will pepper the discussion with examples of singular compactness, including all the classical ones as well as some that are new.  Some examples appear more than once, as cases of different constructions.  The last section of the paper is devoted to a discussion of open questions.

\section{Preliminaries on accessible categories}
\subsection{Accessible functors and ranks of objects}
The next observation is related to the Uniformization Theorem of Makkai and Par\'{e} \cite{MP}, but both the hypotheses and conclusion are stronger.

\begin{prop}  \label{unif}
Let $F:\A\ra\fB$ be a functor preserving filtered colimits, between accessible categories with filtered colimits.  Let $\lambda_1$ be the least rank of accessibility of $\A$ and $\lambda_2$, the least rank of accessibility of $\fB$.  Set
\[  \lambda_3 = \sup\big\{\,\rank(F(X))\;|\; X\in\A\textup{\ is\ }\lambda_1\textup{-presentable} \big\} \, .  \]
Set $\mu_F=\max\{\lambda_1,\lambda_2,\lambda_3\}$.  Then $\mu_F$ is the rank of accessibility of $F$ in the sense of Def.~\ref{func-rank}.
\end{prop}

\begin{proof}
If $\mu$ is such that $\A$ and $\fB$ are $\mu$-accessible, and for all regular $\lambda\geq\mu$, $F$ takes $\lambda$-presentable objects to $\lambda$-presentable ones, then clearly $\mu\geqslant\max\{\lambda_1,\lambda_2,\lambda_3\}$.

Now if a category is $\lambda$-accessible and has filtered colimits, then it is $\mu$-accessible for any $\mu\geq\lambda$, cf. \cite{BR}, 4.1.  Thus $\A$ and $\fB$ are $\mu_F$-accessible.  Let $X\in\A$ be $\lambda$-presentable, $\lambda\geq\mu_F$.  If $\lambda=\lambda_1$ then $\lambda=\mu_F\geq\lambda_3$ so $F(X)$ is certainly $\lambda$-presentable.  If $\lambda>\lambda_1$ then $X$ can be written as a retract of an object that is the colimit of a filtered diagram of $\lambda_1$-presentable objects, along a diagram of size less than $\lambda$.  (See the proof of \cite{BR}, 4.1; $X$ can be written as a $\lambda$-filtered colimit of filtered colimits of $\lambda_1$-presentable objects along diagrams of size less than $\lambda$, thus a retract of the colimit of one such diagram.)  Since $F$ preserves filtered colimits, $F(X)$ is a retract of an object that is the colimit of a filtered diagram of $\lambda_3$-presentable objects, along a diagram of size less than $\lambda$.  Since $\lambda_3\leq\mu_F\leq\lambda$, $F(X)$ is $\lambda$-presentable too.
\end{proof}

\subsection{L\"{o}wenheim-Skolem arguments in presheaf categories}
Let $\C$ be a small category.  The functor category $\Set^\C$ is locally finitely presentable, a fortiori finitely accessible.  For $X\in\Set^\C$, write
\[   u(X) = \coprod_{c\in\C} X(c)  \]
the disjoint union of the sets underlying $X$.

\begin{prop}   \label{presheaf}
For $X\in\Set^\C$, let a subset $S_c\subseteq X(c)$ be given for each object $c\in\C$.  There exists a subfunctor $V\into X$ such that $S_c\subseteq V(c)$ and
\begin{equation}   \label{ineq}
  \card\big(u(V)\big) \leq \max\{\card(\C),\aleph_0\} + \sum_{c\in\C} \card(S_c)
\end{equation}
\end{prop}

\begin{proof}
Let $S_c^{(0)}=S_c$ and having defined $S_c^{(n)}$, set
\[   S_c^{(n+1)} = \big\{y\in X(c)\;|\;
                y=f(x)\textup{\ for some\ }x\in S_{c'}^{(n)}\textup{\ and morphism\ }f:c'\ra c\in\C  \big\}  \]
Note that $S_c^{(n)}\subseteq S_c^{(n+1)}\subseteq X(c)$.  Let $V(c)=\colim_{n\in\N} S_c^{(n)}$ and note that for all morphisms $f:c'\ra c\in\C$ and $x\in V(c')$, $f(x)\in V(c)$.  Thus $V$ is a subfunctor of $X$ and the inequality (\ref{ineq})\ holds.
\end{proof}

Since an object of $\Set^\C$ is finitely presentable if and only if it is a retract of a finite colimit of representables,
\[   \max\{\card(\C),\aleph_0\} = \card(\fp\Set^\C) \, .  \]
By \cite{BR}~Remark~2.3, it follows that for any $X\in\Set^\C$ with $\card(\fp\Set^\C) < \card(u(X))$,
one has $\card(u(X))=\|X\|$.

\begin{cor}  \label{fat}
For $X\in\Set^\C$, let a subset $S_c\subseteq X(c)$ be given for each object $c\in\C$.  Suppose
\[   \card(\fp\Set^\C) < \sum_{c\in\C} \card(S_c)   \]
Then there exists $V\into X$ such that $S_c\subseteq V(c)$ for all $c\in\C$ and
\[      \sum_{c\in\C} \card(S_c) = \|V\| \, .   \]
\end{cor}

\subsection{Lifting smooth chains}
Let $\kappa$ be an ordinal, thought of as an ordered set, hence category.  A \textsl{smooth chain} of shape $\kappa$ in a category $\C$ is a functor $D:\kappa\ra\C$ such that for all limit ordinals $\beta < \kappa$, $D(\beta)$ is a colimit cocone on $D$ restricted to $\{\alpha \;|\; \alpha<\beta \}$.

The next lemma will be used at the end of the proof of the main theorem.  It is the only point in the proof where $F$-structures are extended and highlights the role of that assumption.

\begin{lemma} \label{smooth}
Let $\A$ and $\fB$ be accessible categories with filtered colimits and $F:\A\ra\fB$ a functor preserving filtered colimits.  Let $B:\kappa\ra\fB$ be a smooth chain with colimit $X$.  Suppose that for all $\alpha<\kappa$, the morphism $B(\alpha)\ra B(\alpha+1)$ is in the image of $F$ and $F$-structures extend along morphisms.  Then $X$ is in the image of $F$.
\end{lemma}

\begin{proof}
We will define a smooth chain $A:\kappa\ra\A$ such that $FA$ is isomorphic to $B$.  Use induction on $\alpha<\kappa$.  Select $A(0)\ra A(1)\in\A$ so its $F$-image is isomorphic to $B(0)\ra B(1)$ via $i_0:A(0)\ra B(0)$ and $i_1:A(1)\ra B(1)$.  Given that the diagrams $FA(\beta)$ and $B(\beta)$, $\beta<\alpha$, are isomorphic: for successor $\alpha$, since $F$-structures extend along morphisms, one can select $A(\alpha-1)\ra A(\alpha)$ so its image is isomorphic to $B(\alpha-1)\ra B(\alpha)$ via $i_{\alpha-1}:A(\alpha-1)\ra B(\alpha-1)$ and $i_\alpha: A(\alpha)\ra B(\alpha)$.  At limit $\alpha$, define $A(\alpha)=\colim_{\beta<\alpha} A(\alpha)$ and $i_\alpha$ to be the isomorphism $A(\alpha)\ra B(\alpha)$ induced by the fact that $F$ preserves filtered colimits.  Finally, again since $F$ preserves filtered colimits, $F\big(\colim_{\alpha<\kappa} A(\alpha)\big)$ is isomorphic to $\colim_{\alpha<\kappa} FA(\alpha)$, hence to $X$.
\end{proof}

\subsection{Fixed points in filtered diagrams}
The next few lemmas state that if an object is written as $\kappa$-filtered colimit of $\kappa$-presentable objects in two different ways, and stages of the diagram satisfy a certain continuity property, then stages of the colimit are isomorphic cofinally often.  The proof is a simple back-and-forth argument that can be thought of as a non-additive version of Eklof's lemma~\cite{EM}~Ch.4~Lemma~1.4.

\begin{lemma}  \label{eklof1}
Let $\A$ be a $\lambda$-accessible category with filtered colimits, $U\in\A$ and $\lambda\leq\kappa$ with $\kappa$ uncountable and regular.  Then there exist a diagram $\D$ and functor $K:\D\ra\A$ such that
\begin{itemize}
\item[(1)] $\D$ is a small, $\kappa$-filtered category that has countable filtered colimits (equivalently, colimits of countable chains)
\item[(2)] $K(d)$ is $\kappa$-presentable for all $d\in\D$
\item[(3)] $K$ preserves countable filtered colimits
\item[(4)] the colimit of $K$ is isomorphic to $U$.
\end{itemize}
\end{lemma}

\begin{proof}
Since $\A$ is $\lambda$-accessible with filtered colimits and $\lambda\leq\kappa$, $\A$ is $\kappa$-accessible (\cite{BR}, 4.1).  The canonical diagram of $U$ with respect to $\pres_\kappa(\A)$, the (essentially small)\ category of $\kappa$-presentable objects of $\A$, has the requisite properties.
\end{proof}

\begin{lemma} \label{eklof2}
Let $\kappa$, $\A$ and $U$ be as above, and let $K_i:\D_i\ra\A$, $i=1,2$, be functors satisfying the four conditions of Lemma~\ref{eklof1}.  Let $s_i\in\D_i$, $i=1,2$, be arbitrary.  Then there exist $t_i\in\D_i$, with morphisms $s_i\ra t_i\in\D_i$, such that the colimit inclusions $K_1(t_1)\ra U$ and $K_2(t_2)\ra U$ are isomorphic as objects of $\A/U$.
\end{lemma}

\begin{proof}
We will define chains of objects $s_i^{(0)}\ra s_i^{(1)}\ra\dots\ra s_i^{(n)}\ra\dots$, $n\in\N$, in $\D_i$, such that $s_i^{(0)}=s_i$, for $i=1,2$ respectively, together with a commutative diagram
\[ \xymatrix{ K_1(s_1^{(0)}) \ar[r]\ar[d] & K_1(s_1^{(1)}) \ar[r]\ar[d] & \dots & \dots \ar[r] &
                                            K_1(s_1^{(n)}) \ar[r]\ar[d] & K_1(s_1^{(n+1)}) \ar[d]\ar[r] & \dots \\
K_2(s_2^{(0)}) \ar[r]\ar[ru] & K_2(s_2^{(1)}) \ar[r]\ar[ur] & \dots & \dots \ar[r]\ar[ru] &
                                            K_2(s_2^{(n)}) \ar[r]\ar[ru] & K_2(s_2^{(n+1)}) \ar[r] & \dots
} \]
in $\A$ between their images under $K_i$.  Let $u_i^{(n)}: K_i(s_i^{(n)})\ra U$ denote the colimit structure map; the diagram will also have the property that the composite $K_1(s_1^{(n)})\ra K_2(s_2^{(n)})\llra{u_2^{(n)}}U$ equals $u_1^{(n)}$ and $K_2(s_2^{(n)})\ra K_1(s_1^{(n+1)})\xra{u_1^{(n+1)}}U$ equals $u_2^{(n)}$ for all $n\in\N$.  Set $t_i=\colim_{n\in\N} s_i^{(n)}$.  Note that $K_i(t_i)$ is isomorphic to $\colim_{n\in\N} K(s_i^{(n)})$ and the induced morphism $K_i(t_i)\ra U$ is isomorphic to the colimit of the $u_i^{(n)}$, since $K_i$ preserves filtered colimits.  The morphisms $K_1(t_1)\ra K_2(t_2)$ and $K_2(t_2)\ra K_1(t_1)$ induced by the zig-zags are inverse isomorphisms in $\A/U$.  A morphism $s_i\ra t_i$ exists as the colimit structure map $s_i^{(0)}\ra t_i$, verifying all parts of the claim.

To construct the diagram, having defined $s_1^{(n)}$, the map $u_1(s_1^{(n)})$ factors through $K_2:\D_2\ra\A$ since $K_1(s_1^{(n)})$ is $\kappa$-presentable and $\D_2$ is $\kappa$-filtered.  Without loss of generality, there exists $d_2^{(n)}\in\D_2$ with a morphism $s_2^{(n-1)}\ra d_2^{(n)}$ (when $n>0$)\ and map $K_1(s_1^{(n)})\ra K_2(d_2^{(n)})$ such that the composite $K_1(s_1^{(n)})\ra K_2(d_2^{(n)})\ra U$ equals $u_1(s_1^{(n)})$.

Using the induction hypothesis that $K_2(s_2^{(n)})\ra K_1(s_1^{(n+1)})\xra{u_1^{(n+1)}}U$ equals $u_2^{(n)}$, the maps $K_2(s_2^{(n-1)})\ra K_2(d_2^{(n)})$ and $K_2(s_2^{(n-1)})\ra K_1(s_1^{(n)})\ra K_2(d_2^{(n)})$ become equal after composing with $K_2(d_2^{(n)})\ra U$.  Since $K_2(s_2^{(n-1)})$ is $\kappa$-presentable and $\D_2$ is $\kappa$-filtered, there exists $s_2^{(n)}\in\D_2$ with a morphism $d_2^{(n-1)}\ra s_2^{(n)}$ such that the composites $K_2(s_2^{(n-1)})\ra K_2(d_2^{(n)})\ra K_2(s_2^{(n)})$ and $K_2(s_2^{(n-1)})\ra K_1(s_1^{(n)})\ra K_2(d_2^{(n)})\ra K_2(s_2^{(n)})$ are equal, as required.

Given $s_2^{(n)}$, there is a symmetric argument for the construction of $s_1^{(n+1)}$.
\end{proof}

We will mainly use the following

\begin{cor} \label{eklof3}
Let $F:\A\ra\fB$ be a functor preserving filtered colimits between accessible categories with filtered colimits, with $\mu_F$ as its rank of accessibility, and $\mu_F\leq\kappa$ with $\kappa$ uncountable and regular.  Let $\K_1:\D_1\ra\A$ and $\K_2:\D_2\ra\fB$ be functors satisfying conditions (1)-(2)-(3)\ of Lemma~\ref{eklof1}, with $\colim_{\D_1} K_1 = U$ and $\colim_{\D_2} K_2 = V$, where $F(U)$ is isomorphic to $V$.  Let $s_i\in\D_i$, $i=1,2$, be arbitrary.  Then there exist $t_i\in\D_i$, with morphisms $s_i\ra t_i\in\D_i$, such that the image under $F$ of the colimit inclusion $K_1(t_1)\ra U$ is isomorphic in $\fB/V$ to the colimit inclusion $K_2(t_2)\ra V$.
\end{cor}

\begin{proof}
Apply Lemma~\ref{eklof2} to the composite $FK_1:\D_1\ra\fB$ and $K_2:\D_2\ra\fB$.
\end{proof}

\section{Subobject games}
Let $F:\A\ra\fB$ be a functor with $\fB$ accessible, $X$ any object of $\fB$.  The next definition is motivated by the concept of `strongly almost free module'; see \cite{E}~3.2.

For a regular cardinal $\kappa$, let $\pres_\kappa[X]$ denote the full subcategory of $\Sub(X)$ whose objects are $\kappa$-presentable.  Write $\pres_\kappa[X]\cap F$ for the collection of objects
\[ \big\{ U\in\Sub(X) \;|\; U\textup{\ is\ }\kappa\textup{-presentable and in the image of\ }F \big\}  \]
By abuse of notation, we will write $\pres_\kappa[X]\cap F$ as well for the collection of \textsl{morphisms} in $\pres_\kappa[X]$ that are in the image of $F$.  The context will always make it clear what is meant.
\begin{defn}
A subset $\W$ of $\pres_\kappa[X]\cap F$ \textsl{cannot be blocked} if for all $U\in\W$ and $V\in\pres_\kappa[X]$ there exists $W\in\W$ such that $W$ contains both $U$ and $V$ and $U\into W\in\pres_\kappa[X]\cap F$.
\end{defn}

We wish to prove that under the conditions of the main theorem, $\pres_\kappa[X]$ contains a non-empty set of objects that cannot be blocked.  The proof is indirect and based on Shelah's `subobject game,' as exposed in \cite{H}; see also \cite{E}.

\subsection{The subobject game.} Consider the non-cooperative game of Player $\A$ and Player $\fB$ with rules as follows.  Player $\A$ moves first and chooses an object $A_0\in\pres_\kappa[X]\cap F$.  Player $\fB$ moves second and selects an object $B_1\in\pres_\kappa[X]$ at will.  Thereafter, for each $n\in\N^+$ in turn, Player $\A$ needs to select an object $A_n\in\pres_\kappa[X]\cap F$ that contains both $B_n$ and $A_{n-1}$ and so that $A_{n-1}\into A_n\in\pres_\kappa[X]\cap F$.  In response, Player $\fB$ can select any object $B_{n+1}\in\pres_\kappa[X]$ as he desires.  Player $\A$ loses if she cannot make a valid selection at some $n$ and wins if, against all possible choices by Player $\fB$, the game can continue through all $n\in\N^+$.

A \textsl{game-state before move $n+1$ by Player $\fB$} is thus a chain of objects in $\pres_\kappa[X]$
\[  A_0 \into B_1 \into A_1 \into B_2\into\cdots\into A_{n-1} \into B_n \into A_n   \]
such that the inclusions $A_i\into A_{i+1}$ lie in $\pres_\kappa[X]\cap F$ for each $0\leq i\leq n-1$.

Note that the game is memoryless, i.e.\ the possibilities for Player $\A$'s next move are restricted solely by her previous move and Player $\fB$'s previous move, while the moves of Player $\fB$ are never restricted at all; and the reason for Player $\A$'s losing the game after move $n+1$ is solely her inability to find a suitable $A_{n+1}$, given $A_n$ and $B_n$, regardless of the prior stages of the game.  More precisely, let
\[  \gb = A_0 \into B_1 \into A_1 \into B_2\into\cdots\into A_{n-1} \into B_n \into A_n   \]
and
\[  \gb' = A_0' \into B_1' \into A_1' \into B_2'\into\cdots\into A_{m-1}' \into B_m' \into A_m'   \]
be two game-states, with $A_n$ isomorphic to $A_m'$.  Then any valid continuation of $\gb$, suitably reindexed, is a valid continuation of $\gb'$.  Thus the following property is well-defined:

\begin{defn}
$A\in\pres_\kappa[X]\cap F$ is a \textsl{winning object for Player $\A$} if she can win from any and all game-states whose last object is $A$.  Otherwise, Player $\fB$ can win from any and all game-states whose last object is $A$, and $A$ is said to be a \textsl{losing object for Player $\A$}.
\end{defn}

The next observation is tautologous:

\begin{prop}
Let $\W_\kappa$ be the set of winning objects for Player $\A$.  Then $\W_\kappa$ cannot be blocked in $\pres_\kappa[X]$.
\end{prop}

Indeed, given $U\in\W_\kappa$ and $V\in\pres_\kappa[X]$, consider the set $\cS$ of objects $W\in\pres_\kappa[X]$ such that $W$ contains both $U$ and $V$ and $U\into W\in\pres_\kappa[X]\cap F$.  If $\cS$ were empty, or only contained $W$ that are losing objects for Player $\A$, then $U$ itself would be a losing object for Player $\A$ (with $V$ as winning response by Player $\fB$).  Thus $\cS$ must contain some winning object $W$ for Player $\A$, verifying that the set of winning objects cannot be blocked.

So the onus is to prove $\W_\kappa$ non-empty.  With slight modifications, the proof below would work for any finitely accessible target category $\fB$, and in fact for any category possessing an `underlying set' functor satisfying certain properties, but for the sake of convenience, we prefer to work with presheaf categories at this point.  In particular, we will freely use that $\Sub(X)$ is a complete distributive lattice for any object $X$ of a presheaf category.
The conclusion of the main theorem will be extended to all finitely accessible categories $\fB$ at the end of this section.

\begin{prop}  \label{nonempty}
Let $\A$ be an accessible category with filtered colimits, $\fB$ a functor category $\Set^\C$ where $\C$ is small and $F:\A\ra\fB$ a functor preserving filtered colimits.  Let $X\in\fB$ be an object and $\kappa$ a regular cardinal such that $\max\{\mu_F,\aleph_1\}\leq\kappa$.  Assume that the image of $F$ contains a dense $\kappa$-filter $\F$ of $\Sub(X)$.  Then $\W_\kappa$ is non-empty.
\end{prop}

Excepting the case when the object $X$ itself is assumed to belong to the image of $F$ (which, at least when $\|X\|$ is singular, is the desired conclusion of the singular compactness theorem!)\ it seems hard to give an explicit example of a winning object for Player $\A$.  Prop.~\ref{nonempty} is proved by contradiction.

Assume that $\W_\kappa$ is empty, that is, Player $\A$ has no winning objects.  The initial move $A_0$ of Player $\A$ could be an arbitrary object in $\pres_\kappa[X]\cap F$, and Player $\fB$ must have a winning response to it.  Note that all subsequent moves of Player $\A$ must lie in $\pres_\kappa[X]\cap F$ and since the game is memoryless, a winning response by Player $\fB$ to $A_0$ is also a winning response to the object $A_0$ being played by Player $\A$ at any subsequent state of the game.  Hence, without loss of generality, a winning strategy for Player $\fB$ can be taken to be a map $\s$ from $\pres_\kappa[X]\cap F$ to $\pres_\kappa[X]$ with the property: there exists no sequence of moves $A_n$, $n\in\N$, by Player $\A$ such that
\[   A_0 \into \s(A_0) \into A_1 \into \s(A_1)\into\cdots\into A_n \into \s(A_n)\into\cdots    \]
is a valid game.

Define a smooth chain $C(\alpha)\in \pres_\kappa[X]\cap F$, $\alpha<\kappa$, by induction as follows.  Let $C(0)$ be an arbitrary $\kappa$-presentable object from $\F$.  Thereafter, for successor $\alpha+1$, let $Y_{\alpha+1}$ be the union of $C(\alpha)$ and $\s(C(\alpha))$ in $\Sub(X)$.  By density of $\F$, find a $\kappa$-presentable $C(\alpha+1)\in\F$ containing $Y_{\alpha+1}$.  For limit $\alpha$, let $C(\alpha)$ be the colimit of the chain $C(\beta)$, $\beta<\alpha$.  Since $\F$ is a $\kappa$-filter, $C(\alpha)\in\F$ and is $\kappa$-presentable.

Let $V$ be $\colim_{\alpha<\kappa} C(\alpha)$ in $\Sub(X)$.  Since $\F$ is a $\kappa$-filter, $V$ is isomorphic to $F(U)$ for some $U\in\A$.  Write $U$ as $\colim_{\D_1} K_1$ for some functor $K_1:\D_1\ra\A$ satisfying properties (1)\ through (4)\ of Lemma~\ref{eklof1}.  Cor.~\ref{eklof3} will then apply with $\D_2$ the ordinal (hence poset)\ $\kappa$ and $K_2$ the functor sending $\alpha<\kappa$ to $C(\alpha)$.

We can now define a sequence of moves by Player $\A$ to defeat strategy $\s$.  We will define a countable chain $d_{-1}\ra d_0\ra d_1\ra\cdots\ra d_n\ra\cdots$ of objects and morphisms in $\D_1$, and a countable increasing chain $\alpha_{-1}<\alpha_0<\cdots<\alpha_n<\cdots$ of ordinals less than $\kappa$, by induction.  Set $d_{-1}$ to be an arbitrary object of $\D_1$ and $\alpha_{-1}=0$.  Having defined $d_{n-1}$ and $\alpha_{n-1}$, apply Cor.~\ref{eklof3} to $K_1:\D_1\ra\A$ and $K_2:\kappa\ra\fB$ with $s_1=d_{n-1}$ and $s_2=\alpha_{n-1}+1$.  Set $d_n=t_1$ and $\alpha_n=t_2$ as concluded by the lemma.

Now let $A_n$ be $C(\alpha_n)$ for $n\in\N$.  For all $n\in\N$, $A_{n+1}$ contains $\s(A_n)$ since $C(\alpha_n+1)$ contains $\s(C(\alpha_n))$ and $\alpha_n+1\leq\alpha_{n+1}$.  Finally, since $FK_1(d_n)$ and $K_2(\alpha_n)=C(\alpha_n)=A_n$ are isomorphic as subobjects of $V$, the image of $K_1(d_n)\ra K_1(d_{n+1})\in\A$ by $F$ is isomorphic to the inclusion $A_n\into A_{n+1}$.  Thence Player $\A$'s moves $A_n$, $n\in\N$, defeat the strategy $\s$ of Player $\fB$.  This completes the proof of Prop.~\ref{nonempty}.  \qed

Note that the play by $\A$ that defeats $\s$ is a chain $A_0\into A_1\into A_2\into\cdots$ of subobjects that is the image under $F$ of a composable chain of morphisms in $\A$, namely, $K_1(d_0)\ra K_1(d_1)\ra K_1(d_2)\ra\cdots$ (itself the image under $K_1$ of a composable chain of morphisms in $\D_1$.)  The rules of the game require something less, namely, that each successive inclusion $A_n\into A_{n+1}$ lie in the image of $F$ (adjacent inclusions not necessarily arising as images of morphisms that are composable in $\A$).  This makes the game easier for $\A$ and harder for $\fB$.  Indeed, this is the reason that Player $\fB$'s putative winning strategy needs to depend on the last move by Player $\A$ just as an element of $\Sub(X)$, regardless of which object of $\A$ it is the $F$-image of.  This makes it harder to use the subobject game to build an $\omega$-long chain of subobjects of $X$ arising as the image under $F$ of an $\omega$-long composable chain of morphisms of $\A$.  The assumption that $F$-structures extend along morphisms comes to the rescue, as we will see in the next section.

\section{The big picture}  \label{bigpict}
We return to the hypotheses of Theorem~\ref{main} with the additional assumption that $\fB$ is a functor category $\Set^\C$ with $\C$ small.  We will construct a commutative diagram in $\Sub(X)$:
\[ \xymatrixcolsep{1pc}\xymatrixrowsep{1pc}
\xymatrix{\dots & \dots & \dots && \dots & \dots & \dots & \dots && \dots & \dots \\
B_{0,3} \ar[u]\ar[rr]|(.5)\hole && B_{1,3}\ar[u]\ar[rr] && \dots\ar[r] & B_{\alpha,3}\ar[u]\ar[rr]|(.47)\hole && B_{\alpha+1,3} \ar[u]\ar[rr] &&\dots & B_{\beta,3} \ar[u]\ar[r] & \dots \\
& M_{0,3}\ar[ul]\ar@{-->}[uu]\ar@{-->}[r] & A_{1,3}\ar[u]\ar@/_2pc/@{-->}[uu] &&& A_{\alpha,3} \ar[u] & M_{\alpha,3}\ar[ul]\ar@{-->}[r]\ar@{-->}[uu] & A_{\alpha+1,3} \ar[u]\ar@/_2pc/@{-->}[uu] && \dots && M_{\beta,3}\ar[ul] \\
B_{0,2} \ar[uu]\ar[rr]|(.5)\hole && B_{1,2}\ar[u]\ar[rr]|(.45)\hole && \dots\ar[r] & B_{\alpha,2}\ar[u]\ar[rr]|(.47)\hole && B_{\alpha+1,2} \ar[u]\ar[rr]|(.4)\hole &&\dots & B_{\beta,2} \ar[uu]\ar[r] & \dots \\
& M_{0,2}\ar[ul]\ar@{-->}[uu]\ar@{-->}[r] & A_{1,2}\ar[u]\ar@/_2pc/@{-->}[uu] &&& A_{\alpha,2} \ar[u] & M_{\alpha,2}\ar[ul]\ar@{-->}[r]\ar@{-->}[uu] & A_{\alpha+1,2} \ar[u]\ar@/_2pc/@{-->}[uu] && \dots && M_{\beta,2}\ar[ul] \\
B_{0,1} \ar[uu]\ar[rr]|(.5)\hole && B_{1,1}\ar[u]\ar[rr]|(.45)\hole && \dots\ar[r] & B_{\alpha,1}\ar[u]\ar[rr]|(.47)\hole && B_{\alpha+1,1} \ar[u]\ar[rr]|(.4)\hole &&\dots & B_{\beta,1} \ar[uu]\ar[r] & \dots \\
& M_{0,1}\ar[ul]\ar@{-->}[uu]\ar@{-->}[r] & A_{1,1}\ar[u]\ar@/_2pc/@{-->}[uu] &&& A_{\alpha,1} \ar[u] & M_{\alpha,1}\ar[ul]\ar@{-->}[r]\ar@{-->}[uu] & A_{\alpha+1,1} \ar[u]\ar@/_2pc/@{-->}[uu] && \dots && M_{\beta,1}\ar[ul] \\
B_{0,0} \ar[uu]\ar[rr] && B_{1,0}\ar[u]\ar[rr] && \dots\ar[r] & B_{\alpha,0} \ar[u]\ar[rr] &&          B_{\alpha+1,0} \ar[u]\ar[rr] &&\dots & B_{\beta,0} \ar[uu]\ar[r] & \dots
} \]
The objects $B_{\alpha,n}$ are indexed by the ordinals $\alpha<\cf\|X\|$, $n<\omega$ and form a poset of order type $\cf\|X\|\times\omega$.  Successor columns $\alpha+1$ are preceded by a column of intermediate objects $M_{\alpha,n}$, mapping to $B_{\alpha,n}$ of the predecessor column and also to objects $A_{\alpha+1,n}$ that alternate with the $B_{\alpha+1,n}$.  Limit columns $\beta$ only contain objects $B_{\beta,n}$.  Ordinary arrows indicate inclusions between subobjects of $X$; multiply broken arrows indicate inclusions lying in the image of $F$.  The diagram will have the following properties:
\begin{itemize}
\item[(a)] $\colim_{\alpha<\cf\|X\|} B_{\alpha,0}$ is isomorphic to $X$.
\item[(b)] There is a smooth chain of cardinals $\kappa_\alpha$ such that $\|B_{\alpha,n}\|=\kappa_\alpha$ for all $\alpha<\cf\|X\|$, $n\in\N$.
\item[(c)] There exist objects $a_{\alpha+1,n}\in\A$ for successor $\alpha+1$ and $m_{\alpha,n}\in\A$ for all $\alpha<\cf\|X\|$ and connecting maps $a_{\alpha+1,n-1}\ra a_{\alpha+1,n}$, $m_{\alpha,n-1}\ra m_{\alpha,n}$, $m_{\alpha,n}\ra a_{\alpha+1,n}$ making up a diagram of order type $\{0\ra 1\}\times\omega$ in $\A$ whose image under $F$ is isomorphic to the ladder of morphisms connecting the $M_{\alpha,n}$ to the $A_{\alpha+1,n}$, for fixed $\alpha$ and $n\in\N$.
\item[(d)] $m_{\alpha,n}$ is $\kappa_\alpha^+$-presentable for $\alpha<\cf\|X\|$, $n\in\N^+$.
\item[(e)] $M_{\alpha,n}$ is the intersection of $A_{\alpha+1,n}$ and $B_{\alpha,n}$ in $\Sub(X)$, for all $\alpha<\cf\|X\|$, $n\in\N^+$.
\item[(f)] For all limit ordinals $\beta<\cf\|X\|$ and $n\in\N^+$, $B_{\beta,n-1}$ is contained in $\colim_{\alpha < \beta} B_{\alpha,n}$.
\end{itemize}

\noindent
Let us show that the conclusion that $X$ belongs to the image of $F$ follows from (a),(c),(e),(f).  Taking the colimit in $\Sub(X)$ along columns of the big picture, one obtains a diagram
\[  \xymatrixcolsep{1pc}\xymatrixrowsep{1pc}\xymatrix{
B(0) \ar[rr] && B(1) \ar[rr] && B(2) \ar[r] & \dots \ar[r] & B(\alpha) \ar[rr] && B(\alpha+1) \ar[r] & \\
& M(0) \ar[ul]\ar@{-->}[ur] && M(1) \ar[ul]\ar@{-->}[ur] && \dots && M(\alpha) \ar[ul]\ar@{-->}[ur]
} \]
of subobjects of $X$, where $B(\alpha)=\colim_n B_{\alpha,n}$ and $M(\alpha)=\colim_n M_{\alpha,n}$.  The map $B(\alpha)\la M(\alpha)$ is an isomorphism.  Indeed, $B_{\alpha,n}$ is a subobject both of $B_{\alpha,n+1}$ and $A_{\alpha+1,n+1}$ hence, by (e), of their intersection $M_{\alpha,n+1}$.  The inclusions $B_{\alpha,n}\ra M_{\alpha,n+1}$ induce the map $B(\alpha)\ra M(\alpha)$ inverse to $B(\alpha)\la M(\alpha)$.  Each connecting map
\[   M(\alpha)\ra B(\alpha+1) = \colim_n B_{\alpha+1,n} = \colim_n A_{\alpha+1,n}    \]
is in the image of $F$ by (c).  For all limit ordinals $\beta<\cf\|X\|$, $\colim_{\alpha<\beta} B(\alpha)$ is isomorphic to $B(\beta)$ by (f).  Lemma~\ref{smooth} applies now to show that $\colim_{\alpha<\cf\|X\|} B(\alpha)$ is in the image of $F$.  But that colimit is isomorphic to $X$, since it is a subobject of $X$ and (a)\ holds. \qed

Let us construct the diagram now.  To define the bottom row, we resume the notation of Prop.~\ref{presheaf}.  Since $\card(\fp\Set^\C)<\|X\|$, one has $\card(u(X)) = \|X\|$.  Write $u(X)$ as the union of a smooth chain of subsets $S_\alpha$, $\alpha<\cf\|X\|$, with $\card(S_\alpha)<\|X\|$ for all $\alpha$.  Without loss of generality,
\[   \max\{\cf\|X\|,\mu_F,\card(\fp\Set^\C),\phi_{F,X}\} < \card(S_0)   \]
where $\mu_F$ is the rank of accessibility of $F$ and $\phi_{F,X}$ the threshold cardinal from assumption (ii)\ of Theorem~\ref{main}.  Using Cor.~\ref{fat}, define the objects $B_{\alpha,0}$ by transfinite induction:
\begin{itemize}
\item let $B_{0,0}$ be a subobject of $X$ such that $S_0\subseteq u(B_{0,0})$ and $\card(S_0)=\|B_{0,0}\|$
\item for successor $\alpha$, let $V_\alpha$ be a subobject of $X$ such that $S_\alpha\subseteq u(V_\alpha)$ and $\card(S_\alpha)=\|V_\alpha\|$, then let $B_{\alpha,0}$ be the union of $B_{\alpha-1,0}$ and $V_\alpha$ in $\Sub(X)$
\item for limit ordinals $\beta$, let $B_{\beta,0}$ be $\colim_{\alpha<\beta} B_{\alpha,0}$ in $\Sub(X)$.
\end{itemize}
(a)\ is therefore satisfied.  Write $\kappa_\alpha=\|B_{\alpha,0}\|$; this is a smooth chain of cardinals for $\alpha<\cf\|X\|$.  For each successor ordinal $\alpha+1<\cf\|X\|$ choose a set $\W_{\alpha+1}$ of objects that cannot be blocked in $\pres_{\kappa_{\alpha+1}}[X]$.

The diagram is now constructed by induction on $n$.  Whenever the objects $B_{\alpha,n}$ have been defined for some $n\in\N$, choose, for each limit ordinal $\beta<\cf\|X\|$, an arbitrary enumeration $\{ b_{\beta,n}(i)\;|\; i < \kappa_\beta\}$ of the elements of $u(B_{\beta,n})$.

To find $A_{\alpha+1,1}$ for a successor ordinal $\alpha+1$: since $\W_{\alpha+1}$ cannot be blocked in $\pres_{\kappa_{\alpha+1}}[X]$, there exists $A_{\alpha+1,1}\in\W_{\alpha+1}$ containing $B_{\alpha+1,0}$.  Choose also some $a_{\alpha+1,1}\in\A$ and isomorphism $i_{\alpha+1,1}: F(a_{\alpha+1,1}) \ra A_{\alpha+1,1}$.

To find $A_{\alpha+1,n}$ for $n>1$: again since $\W_{\alpha+1}$ cannot be blocked in $\pres_{\kappa_{\alpha+1}}[X]$, there exists $A_{\alpha+1,n}\in\W_{\alpha+1}$ containing $B_{\alpha+1,n-1}$ (a fortiori $A_{\alpha+1,n-1}$)\ and such that $A_{\alpha+1,n-1}\into A_{\alpha+1,n}\in \pres_{\kappa_{\alpha+1}}[X]\cap F$.  Since $F$-structures are assumed to extend along morphisms, there exists $a_{\alpha+1,n-1}\ra a_{\alpha+1,n}\in\A$ and isomorphism $i_{\alpha+1,n}: F(a_{\alpha+1,n}) \ra A_{\alpha+1,n}$
making
\[ \xymatrix{  F(a_{\alpha+1,n-1}) \ar[r]\ar[d]^{i_{\alpha+1,n-1}} &  F(a_{\alpha+1,n})\ar[d]^{i_{\alpha+1,n}} \\
A_{\alpha+1,n-1} \ar@{>->}[r] & A_{\alpha+1,n}
} \]
commute.

Express now $a_{\alpha+1,n}$ as $\colim_{\D_1} K_1$ for a functor $K_1:\D_1\ra\A$ satisfying the hypotheses of Lemma~\ref{eklof1}, with $\kappa_\alpha$ playing the role of the $\kappa$ of Lemma~\ref{eklof1}.  Since $m_{\alpha,n-1}$ is $\kappa_\alpha^+$-presentable, the composite $m_{\alpha,n-1} \ra a_{\alpha+1,n-1} \ra a_{\alpha+1,n}$ factors through $\D_1$, say, through $s_1\in\D_1$.  (When $n=1$, choose an arbitrary object $s_1\in\D_1$.)

For $\alpha<\cf\|X\|$, $n\in\N^+$, use Cor.~\ref{fat} to find a subobject $U_{\alpha,n}$ of $X$ containing
\[\begin{cases}
\{ b_{\beta,n-1}(i) \}\textup{\ for all limit ordinals\ }\beta\textup{\ with\ }\alpha < \beta < \cf|X|,\textup{\ all\ } i <\kappa_\alpha \\
B_{i,n}\textup{\ for all\ }i < \alpha \\
A_{\alpha,n-1}\textup{\ when\ }\alpha\textup{\ is a successor ordinal} \\
B_{\alpha,n-1}\textup{\ when\ }\alpha\textup{\ is a limit ordinal}.
\end{cases}\]
such that $\|U_{\alpha,n}\|=\kappa_\alpha$.  Let $A_{\alpha+1,n} \cap U_{\alpha,n}$ denote the intersection of $A_{\alpha+1,n}$ and $U_{\alpha,n}$ in $\Sub(X)$.  Express $A_{\alpha+1,n}$ as $\colim_{\D_2} K_2$, the $\kappa_\alpha^+$-directed union of its $\kappa_\alpha^+$-presentable subobjects.  Note that $A_{\alpha+1,n} \cap U_{\alpha,n}$ is such a subobject, corresponding, say, to $s_2\in\D_2$.

Cor.~\ref{eklof3} applies now to give morphisms $s_1\ra t_1\in\D_1$ and $s_2\ra t_2\in\D_2$ such that $F\big(K_1(t_1)\ra a_{\alpha+1,n}\big)$ and $K_2(t_2)\ra A_{\alpha+1,n}$ are isomorphic as subobjects of $A_{\alpha+1,n}$.  Set $m_{\alpha,n}$ to be $K_1(t_1)$, so (d)\ holds.  The connecting map $m_{\alpha,n} \ra a_{\alpha+1,n}$ is the colimit inclusion and (for $n>1$)\ the map $m_{\alpha,n-1} \ra m_{\alpha,n}$ is given by the composite $m_{\alpha,n} \ra K_1(s_1) \ra K(t_1)$.  Let $M_{\alpha,n}$ be $F(m_{\alpha,n})$.  This verifies condition (c).

Define $B_{\alpha,n}$ to be the union of $U_{\alpha,n}$ and $M_{\alpha,n}$ in $\Sub(X)$.  Note that $B_{\alpha,n}$ includes $B_{i,n}$ for $i<\alpha$ and $A_{\alpha,n-1}$ (resp.\ $B_{\alpha,n-1}$), hence fits into the diagram.  $\|B_{\alpha,n}\|=\kappa_\alpha$, so (b)\ is maintained.  Since $M_{\alpha,n}$ is isomorphic to $K_2(t_1)$, it is a subobject of $A_{\alpha+1,n}$ containing $A_{\alpha+1,n} \cap U_{\alpha,n}$.  Since in a presheaf category the lattice of subobjects of any object is distributive,
\[ A_{\alpha+1,n} \cap B_{\alpha,n} = A_{\alpha+1,n} \cap (U_{\alpha,n} \cup M_{\alpha,n}) =
   (A_{\alpha+1,n} \cap U_{\alpha,n} ) \cup (A_{\alpha+1,n} \cap M_{\alpha,n}) = M_{\alpha,n} \]
Thence condition (e)\ is satisfied.

Finally, $B_{\alpha,n}$ contains the elements $b_{\beta,n-1}(i)$ for limit ordinals $\beta$ with $\alpha < \beta < \cf\|X\|$, $i <\kappa_\beta$.  Since for limit $\beta$, one has $\kappa_\beta=\sup\{\kappa_\alpha\;|\;\alpha<\beta\}$, (f)\ holds.   \qed

There remains the task of removing the assumption, stated at the beginning of section~\ref{bigpict}, that the target category $\fB$ in the main theorem is a functor category $\Set^\C$.  First, a definition from~\cite{BR} and an easy observation.

\begin{defn}
A functor $G:\A\ra\fB$ is \textsl{iso-full} if for every isomorphism $i\in\fB$ there exists an isomorphism $i_0\in\A$ such that $G(i_0)=i$.
\end{defn}

\begin{prop}     \label{composite}
$(a)$ If $G$ is iso-full then $G$-structures extend along morphisms.  $(b)$ If $F:\A\ra\fB$ and $G:\fB\ra\C$ are functors such that $F$-structures extend along morphisms and $G$ is iso-full, then $GF$-structures extend along morphisms.
\end{prop}

\begin{proof}
$(a)$ Given $g:X\ra Y$ and an isomorphism $i:G(U)\ra G(X)$, let $i_0$ be a morphism such that $G(i_0)=i$ and set $f=gi_0$ and $j=\id_{G(Y)}$ to satisfy Def.~\ref{extend}.  (Note that this part does not require that $i_0$ be an isomorphism.)  $(b)$ Given $g:X\ra Y$ and isomorphism $i:GF(U)\ra GF(X)$, find an isomorphism $i_0:F(U)\ra F(X)$ such that $G(i_0)=i$, then $f:U\ra V$ and isomorphism $j:F(V)\ra F(Y)$ such that the diagram commutes.
\end{proof}

The next lemma gives sufficient conditions on $\fB$ for the conclusion of the main theorem to stand.  These conditions state, essentially, that subobject structures of $\fB$ become similar to subobject lattices of presheaf categories, via a comparison functor $G:\fB\ra\Set^\C$.  The most technical-looking condition, Lemma~\ref{general}\,(6), is that an analogue of the downward L\"{o}wenheim--Skolem theorem holds.  Conditions (5)\ and (6)\ turn out to be redundant, i.e.\ to be consequences of slight variants of conditions (1)\ through (4).  We prefer to postpone discussing this and to use the lemma to quickly and directly prove our main theorem first.

\begin{lemma}   \label{general}
Let $\fB$ be an accessible category with filtered colimits such that there exists a small category $\C$, a functor $G:\fB\ra\Set^\C$ with rank of accessibility $\mu_G$ and cardinal $\lambda_G\geq\mu_G$ with the following properties:
\begin{itemize}
\item[(1)] $G$ is iso-full
\item[(2)] $G$ reflects split epimorphisms
\item[(3)] $G$ preserves filtered colimits
\item[(4)] $G$ preserves monomorphisms
\item[(5)] for all objects $X$ of $\fB$, the functors $\Sub(X)\ra\Sub(G(X))$ induced by $G$ are full (that is, if $U$ and $V$ are subobjects of $X$ such that $G(U)$ in included in $G(V)$ then $U$ is included in $V$)
\item[(6)] for all regular cardinals $\kappa\geq\lambda_G$, object $X\in\fB$ and $\kappa$-presentable subobject $V\into G(X)$, there exists a $\kappa$-presentable subobject $U\into X$ such that $G(U)$ contains $V$.
\end{itemize}

Let now $\A$ be an accessible category with filtered colimits, $F:\A\ra\fB$ a functor preserving filtered colimits with rank of accessibility $\mu_F$ and object $X\in\fB$ such that
\[  \max\big\{\lambda_G,\card(\C),\mu_F\big\} < \|X\|  \]
and assumptions (i),(ii),(iii)\ of Theorem~\ref{main} hold.  Then $X$ is in the image of $F$.
\end{lemma}

\begin{proof}
Consider the composite $GF:\A\ra\fB\ra\Set^\C$.  We want to apply the already proved case of Theorem~\ref{main} to $GF$ and the object $G(X)$.  Since $\mu_{GF}\leq\max\{\mu_F,\mu_G\}$, we have
\[   \max\{\mu_{GF},\card(\fp\Set^\C)\}\leq\max\{\mu_F,\card(\C)\} < \|X\| \, . \]
Since $G$ preserves filtered colimits and reflects split epimorphisms, Propositions~4.3 and 3.7 of \cite{BR} imply that $\|G(U)\|=\|U\|$ for all objects $U$ with $\|U\|\geq\mu_G$ and thus $\|G(X)\|=\|X\|$ is singular.

Let now $\kappa$ be a regular cardinal greater than or equal to $\mu_G$, $\F$ a $\kappa$-filter of $\Sub(X)$ and $G(\F)$ its image in $\Sub(G(X))$.  Consider any well-ordered chain $C$ of $\kappa$-presentable objects in $G(\F)$ of length less than or equal to $\kappa$.  Thanks to assumption (5), there exists a pre-image of this chain in $\F$ consisting of $\kappa$-presentable objects.  Since $\F$ is a $\kappa$-filter and $G$ preserves filtered colimits, the colimit of $C$ exists in $G(\F)$.  This implies that $G(\F)$ is a $\kappa$-filter too.

Assumption (6)\ implies that $G$ takes $\kappa$-dense filters to $\kappa$-dense filters for all regular $\kappa\geq\lambda_G$.  By the assumption on $F$, there exists a cardinal $\phi_{F,X}<\|X\|$ such that for all regular $\kappa\geq\phi_{F,X}$, the image of $F$ contains a dense $\kappa$-filter of $\Sub(X)$.  Set $\phi_{GF,G(X)}=\max\{\lambda_G,\phi_{F,X}\}<\|G(X)\|$.  Then assumption (ii)\ of Theorem~\ref{main} is satisfied for $GF$ too.  $GF$-structures extend along morphisms by Prop.~\ref{composite}.  The known case of Theorem~\ref{main} thus applies to $GF:\A\ra\fB\ra\Set^\C$ to show that $G(X)$ lies in the image of $GF$.  Since $G$ is iso-full, this implies that $X$ lies in the image of $F$.
\end{proof}

\begin{prop}  \label{main1}
Any finitely accessible category $\fB$ satisfies the assumptions of Lemma~\ref{general} with $\C=(\fp\fB)^\op$ and $\lambda_G=\card(\fp\fB)$.
\end{prop}

\begin{proof}
Consider the restricted Yoneda embedding $G:\fB\ra\Set^\C$, sending $X\in\fB$ to the presheaf $\hom_\fB(-,X)$.  Since $G$ is a full embedding, conditions (1), (2)\ and (5)\ hold.  (3)\ is well-known and (4)\ holds for all Yoneda embeddings.  Note that $\mu_G=\aleph_0$ by Prop.~\ref{unif} and $G$ preserves the ranks of all objects with uncountable rank.  Also note that $\fB$ is equivalent to the category of flat presheaves on $\fp\fB$.  This theory can be axiomatized by $\card(\fp\fB)$ many sentences of the logic $\LL_{\card(\fp\fB)^+,\,\omega_0}$, cf.~\cite{MP} 4.4.3 and 3.2.3, or~\cite{presh} (2.1)\ through (2.5)\ for an explicit set of formulas.  Condition (6) now follows as in the proof of the downwards L\"{o}wenheim-Skolem theorem; see e.g.\ \cite{MP}~3.3.1.

This completes the proof of Theorem~\ref{main}.
\end{proof}

We do not have a direct characterization of the class of categories that can play the role of the target $\fB$ in Theorem~\ref{main}, but it is properly larger than the class of finitely accessible categories.  Recall that to any Abstract Elementary Class in Shelah's sense \cite{AEC}, one can associate the category whose objects are models and whose morphisms are strong embeddings.  Let us refer to these as AEC categories.  Theorem 5.5 of \cite{BR} implies that any AEC category $\fB$ satisfies the assumptions of Lemma~\ref{general} via the inclusions $\fB\ra\emb(\Sigma)\ra\Set^\C$ where $\Sigma$ is the signature of the AEC, $\emb(\Sigma)$ is the finitely accessible category of $\Sigma$-structures and their embeddings, and $\C$ is $\big(\fp\emb(\Sigma)\big)^\op$.  (The existence of a L\"{o}wenheim--Skolem cardinal is one of the axioms of Abstract Elementary Classes.)  Not all AEC categories are finitely accessible; in fact, the least rank of accessibility of an AEC category can be as large as desired.

For future reference, let us include Lemma~\ref{main2}, an improved version of Lemma~\ref{general}.  It involves a smaller set of conditions but its proof relies on some advanced machinery from the 2-category of accessible categories.  It will not be needed in the rest of this paper; in all our applications here, the category $\fB$ will be finitely accessible (and often, locally finitely presentable).

\begin{lemma}   \label{main2}
Let $\fB$ be an accessible category with filtered colimits such that there exists a small category $\C$ and functor $G:\fB\ra\Set^\C$ with the following properties:
\begin{itemize}
\item[(1)] $G$ is iso-full
\item[($2^+$)] $G$ is nearly full (see \cite{BR}~Def.~5.1 for this notion)
\item[(3)] $G$ preserves filtered colimits
\item[($4^+$)] $G$ preserves monomorphisms and is faithful.
\end{itemize}

These data satisfy the conditions of Lemma~\ref{general} with $\lambda_G=\max\{\card(C)^+,\mu_G\}$.
\end{lemma}

\begin{proof}
A nearly full and faithful functor reflects split epimorphisms; see \cite{BR}~Remark~5.2.  It is immediate that a nearly full functor preserving monomorphisms satisfies (5).

To verify the L\"{o}wenheim-Skolem property (6), consider the pullback of categories
\[ \xymatrix{
\fB_0 \ar[r]\ar[d]^{G_0} & \fB \ar[d]^G \\
(\Set^\C)_\mono \ar[r]^i & \Set^\C
} \]
where $(\Set^\C)_\mono$ is the category with the same objects as $\Set^\C$ but morphisms the monomorphisms, and $i$ the inclusion.  Said more directly, $\fB_0$ can be considered as the subcategory of $\fB$ with the same objects as $\fB$, its morphisms being those morphisms of $\fB$ that are sent to monomorphisms by $G$.

Since faithful functors reflect monomorphisms, all morphisms of $\fB_0$ are monomorphisms as well.  Since $\fB$ has filtered colimits and $G$ preserves them, and filtered colimits of monomorphisms are monos in $\Set^\C$, $\fB_0$ has filtered colimits and they are preserved by $G_0$.

Since $i$ is a transportable functor, this pullback square is equivalent to a pseudopullback.  (See \cite{MP}~Prop.~5.1.1.)  By the Limit Theorem of Makkai and Par\'{e}, $\fB_0$ is accessible and \cite{CR} gives a bound for the degree of accessibility of $\fB_0$.

Indeed, $(\Set^\C)_\mono$ and $\Set^\C$ are finitely accessible and $i$ preserves filtered colimits.  Let $u(X)$ denote the disjoint union of the sets underlying $X$, for a functor $X:\C\ra\Set$.  Note that both in $(\Set^\C)_\mono$ and $\Set^\C$, one has $\|X\|=\card(u(X))^+$ whenever $\card(\C)<\card(u(X))$.  Hence $i$ preserves $\kappa$-presentable objects for $\kappa>\card(\C)$.  By the proof of Theorem 3.1 of \cite{CR}, $\fB_0$ is $\lambda_G=\max\{\card(C)^+,\mu_G\}$-accessible.  Since it has filtered colimits, it is well-$\lambda_G$-accessible.

Let a regular $\kappa\geq\lambda_G$ be given now, with an object $X\in\fB$ and $\kappa$-presentable subobject $V\into G(X)$.  Considering $X$ as an object of $\fB_0$, write it as the $\kappa$-directed supremum of its $\kappa$-presentable subobjects.  Write $G(X)$ as the $\kappa$-directed supremum of its $\kappa$-presentable subobjects containing $V$.  Applying the fixed-point Lemma~\ref{eklof3} to $G$, there exists a $\kappa$-presentable $U\into X$ such that $G(U)$ contains $V$, as desired.
\end{proof}

For finitely accessible $\fB$ with $G$ the restricted Yoneda embedding, Lemma~\ref{main2} specializes to Prop.~\ref{main1}.  Note that the existence of a functor into a presheaf category satisfying conditions (1), ($2^+$) and (3)\ of Lemma~\ref{main2} characterizes what we called \textsl{abstract elementary categories} in~\cite{BR}.  These seem to share all key structural properties of Abstract Elementary Classes, but their morphisms are not restricted to be monomorphisms.  Conversely, Remark~5.8\,(2)~of~\cite{BR} shows how to associate an AEC to an abstract elementary category.

\subsection{Comparison with the work of Shelah and Hodges.}  We discuss how Theorem~\ref{main} captures the original examples of singular compactness, and add some new ones.  First, an important family of functors $F:\A\ra\fB$ so that $F$-structures extend along morphisms.  We say that \textsl{subobjects have complements} in a category $\A$ if coproduct inclusions are monomorphisms and conversely, for all monomorphisms $X\into Y$, there exists a monomorphism $Z\into Y$ such that $X\ra Y\la Z$ is a coproduct.

\begin{lemma}  \label{boolean}
Suppose subobjects have complements in $\A$ and the functor $F:\A\ra\fB$ preserves finite coproducts.  Let $\A_\mono$ be the subcategory of $\A$ with the same objects as $\A$, but with morphisms the monomorphisms.  Then $F$-structures extend along morphisms for the functor $F:\A_\mono\ra\fB$.
\end{lemma}

\begin{proof}
Given $g:X\into Y$ and $U$ of $\A$, together with an isomorphism $i:F(U)\ra F(X)$, let $Z\into Y$ be a complement of $g:X\into Y$ and define $V$ to be the coproduct of $U$ and $Z$, with $f:U\into V$ the coproduct inclusion.  Consider the diagram
\[ \xymatrix{  F(U) \ar[r]\ar[d]^i &  F(U)\sqcup F(Z) \ar[d]^{i\sqcup\id}\ar[r]^(.65)v & F(V) \\
 F(X) \ar[r] & F(X)\sqcup F(Z) \ar[r]^(.65)y & F(Y) } \]
The isomorphism $j:F(V)\ra F(Y)$ needed to make the diagram of Def.~\ref{extend} commute is $y\circ(i\sqcup\id)\circ v^{-1}$.
\end{proof}

Examples of categories with complemented subobjects are $\Set$, the functor category $\Set^G$ where $G$ is a groupoid, $\sh(\B)$, the category of sheaves on a complete boolean algebra $\B$ equipped with its canonical topology, and the category of $k$-vector spaces for a field $k$ (or more generally, of $k$-modules for a division ring $k$).

\begin{example}
For any one-sorted, equational variety of universal algebras, let $\Alg$ be the category of algebras and homomorphisms and $F:\Set\ra\Alg$ the free algebra functor.  Theorem~\ref{main} applies to the restriction $F:\Set_\mono\ra\Alg$; the threshold cardinal $\max\{\mu_F,\card(\fp\Alg)\}$ is the greater of $\aleph_0$ and the number of function symbols in the algebra.

Setting $\Alg$ to be groups or abelian groups, this recovers the historically first cases of singular compactness.  As Hodges~\cite{H} remarks, in this context one should look for applications of singular compactness only among varieties with the Schreier property, that is, varieties with the property that any subalgebra of a free algebra is free.  A well-known theorem of Neumann et al.~\cite{sch1}~\cite{sch2} asserts that the only varieties of groups with the Schreier property are groups, abelian groups, and abelian groups of exponent $p$, for a prime $p$.
\end{example}

\begin{example}
Let $k$ be a field, $\Vect_k$ the category of $k$-vector spaces, $\Alg_k$ a category of $k$-linear algebras, and $F:\Vect_k\ra\Alg_k$ the free algebra functor again.  Theorem~\ref{main} applies to $F:(\Vect_k)_\mono\ra\Alg_k$.  No classification of Schreier varieties is known, but classical theorems of Kurosh, Shirshov and Witt give several examples: absolutely free algebras (also known as `non-associative algebras,' i.e.\ algebras with a set of binary, $k$-linear operations and no identities), commutative (but not associative)\ algebras, anti-commutative algebras, and Lie algebras and Lie $p$-algebras.  See~\cite{sch} for an excellent overview.

The corresponding examples of singular compactness seem to be new.  We do not know, however, whether the statement is sharp here; for example, whether there exists a non-free Lie algebra whose size is a regular cardinal, all of whose sub-Lie-algebras of lesser, regular cardinality are free.  (Such examples are known to exist for groups and abelian groups, cf.~\cite{EM}.)

Almost free modules over a ring $R$ correspond to the case of the free functor $F:\Set_\mono\ra\Mod_R$.  Note that for an object $X$ of $\Mod_R$, its size $\|X\|$ is the least cardinal $\kappa$ such that $X$ is $\kappa$-generated as an $R$-module, just as in \cite{EM}.
\end{example}

\begin{example}   \label{decomp}
Let $\fB$ be a locally finitely presentable category and $U_i$, for $i\in I$, a set of objects of $\fB$.  Let the functor $F:\Set^I\ra\fB$ take the collection of sets $\{X_i\}_{i\in I}$ to $\bigsqcup_{i\in I}\sqcup_{X_i} U_i$, that is, the coproduct of the $X_i$-fold copowers of the $U_i$.  Theorem~\ref{main} applies to the restriction $F:(\Set^I)_\mono\ra\fB$, with threshold cardinal the greater of the supremum of $\rank(X_i)$, $i\in I$, and $\card\fp\fB$.  $X$ belongs to the image of $F$ if and only if it is $U_i$-decomposable, that is, isomorphic to the coproduct of a set of objects, each of which is one of the $U_i$.  Setting $\fB$ to be the category of $R$-modules, one obtains the usual notion of decomposable modules.
\end{example}

The following example is also due to Hodges~\cite{H}, who treats it in a slightly different form.  Here the target category is finitely accessible, but not locally finitely presentable.

\begin{example}
Let $k$ be any field, let $\Field_k$ be the category of field extensions of $k$ and homomorphisms, and let $F:\Set_\mono\ra\Field_k$ take the set $I$ to the purely transcendental extension of $k$ generated by the set of indeterminates $x_i$ for $i\in I$.  Note that $F$ is left adjoint to the functor that sends a field $L\in\Field_k$ to the set $L\setminus k$.  In particular --- while not all coproducts exist in $\Field_k$ --- coproducts in $\Set$ are preserved by $F$, and Lemma~\ref{boolean} applies.

A field $L$ is in the image of $F$, that is, isomorphic to $F(I)$ for some set $I$, if and only if it is a purely transcendental extension of $k$.  For an object $X/k$ of $\Field_k$, note that $\|X/k\|$ is the least cardinal $\kappa$ such that $X$ is $\kappa$-generated as an extension of $k$.  Singular compactness thus says in this case: if $L$ is an overfield of $k$ that is singularly generated over $k$, and such that all intermediate fields between $k$ and $L$ that are regularly generated over $k$ are purely transcendental, then $L$ itself is purely transcendental.

It is not clear to us whether this theorem has an applicable case.  L\"{u}roth's theorem states that all fields intermediate between $k$ and $k(x)$, for a single indeterminate $x$, are purely transcendental extensions of $k$ of degree 1, and the theorem of Castelnuovo-Zariski states that when $k$ is algebraically closed, all fields intermediate between $k$ and $k(x,y)$, for indeterminates $x$ and $y$, are purely transcendental extensions of $k$.  (This is not necessarily so for $k$ not algebraically closed.)
\end{example}

Singular compactness is sometimes said to be about ``an abstract notion of free''. (Cf.\ the first sentence of Eklof~\cite{E}, or the last section of Hodges~\cite{H}.)  This is indeed so in the above examples, where the functor $F$ is a left adjoint, i.e.\ free construction.  But there are cases of the main theorem where the intuition of ``free structure'' is absent.

\begin{example}     \label{lingraph}
Let $\lingraph$ be the category whose objects are triples $\lc X,R,\prec\rc$ where $X$ is a set and $R$ a symmetric relation on $X$ (i.e.\ $(X,R)$ is a graph)\ and $\prec$ is a linear order on the set of vertices $X$.  A morphism in $\lingraph$ is to be a monotone, downward closed embedding of subgraphs.  That is, a morphism from $\lc X,R,\prec_X\rc$ to $\lc Y,Q,\prec_Y\rc$ is an injective map $g:X\ra Y$ such that $uRv$ if and only if $g(u)Qg(v)$; if $u \prec_X v$ then $g(u) \prec_Y g(v)$; and if $v$ is in the image of $g$ and $u\prec_Yv$, then $u$ is in the image of $g$.  Let $\graph$ be the category of graphs and embeddings, and $F:\lingraph\ra\graph$ the functor that forgets the ordering of vertices.  It is immediate that $F$-structures extend along morphisms.  Indeed, let $g$ be a morphism from $\lc X,R,\prec_X\rc$ to $\lc Y,Q,\prec_Y\rc$ and let $i:(U,S)\ra(X,R)$ be an isomorphism of graphs.  Let $V$ be the disjoint union of $U$ and $Y\setminus g(X)$ and $j:V\ra Y$ the bijection that is $g\circ i$ on $U$ and the identity on $Y\setminus g(X)$.  Define an order on $V$ by setting $r\prec_V s$
\[ \begin{cases}
\textup{for all\ }r,s\in U\textup{\ such that\ }r\prec_U s \\
\textup{for all\ }r\in U\textup{\ and\ }s\in Y\setminus g(X) \\
\textup{for all\ }r,s\in Y\setminus g(X)\textup{\ such that\ }r\prec_Y s.
\end{cases} \]
The graph with linearly ordered vertices $\lc V,j^{-1}Q,\prec_V\rc$ and $f:\lc U,S,\prec_U\rc\ra\lc V,j^{-1}Q,\prec_V\rc$ which is the identity on $U$, complete the diagram of Def.~\ref{extend}.

Fix a regular uncountable cardinal $\mu$, and consider the full subcategory $\mu$-$\lingraph$ of $\lingraph$ whose objects $\lc X,R,\prec\rc$ are also supposed to satisfy that for all vertices $v\in X$,
\begin{equation}   \label{mu-short}
   \card\{u\in X\;|\;u\prec v\textup{\ and\ }uRv\} < \mu  \; .
\end{equation}
$F$-structures extend for the forgetful functor $F:\mu$-$\lingraph\ra\graph$, as shown by the same construction.

The category $\mu$-$\lingraph$ has directed colimits, computed as directed unions on underlying graphs.  (The condition that morphisms be downward-closed monotone embeddings ensures that (\ref{mu-short})\ is satisfied for the directed union too.)  Thanks to (\ref{mu-short}), any object of $\mu$-$\lingraph$ is the $\mu$-directed colimit of its restrictions to downward closed subgraphs of size less than $\mu$.  Since there is only a set of isomorphism classes of objects in $\mu$-$\lingraph$ with underlying set of size less than $\mu$, and these objects are $\mu$-presentable, $\mu$-$\lingraph$ is a $\mu$-accessible category with filtered colimits.  $\graph$ is finitely accessible and $F$ preserves filtered colimits; thus Theorem~\ref{main} applies.
\end{example}

This method of rigidifying models by adding an underlying linear order is not limited to graphs, but works more generally in the category $\str(\Sigma)$ of structures with signature $\Sigma$.  Since any Abstract Elementary Class is naturally a subcategory of some $\str(\Sigma)$, one can obtain examples of functorial singular cardinal compactness with values in AEC's this way.

\begin{example}    \label{wellgraph}
Let $\mu$-$\wellgraph$ be the subcategory of $\mu$-$\lingraph$ where the ordering of vertices is required to be a well-order.  Keep the rest of the set-up of Example~\ref{lingraph}.  The same arguments work, with small changes.  Shelah says that such a graph has coloring number $\mu$; see also Hodges~\cite{H}.
\end{example}

Transfinite induction shows that any object of $\mu$-$\wellgraph$ is $\mu$-colorable.  What of graph colorability without well-ordered vertices?  Let $\mu$-$\graph$ be the category of $\mu$-colored graphs and embeddings (preserving the coloring) and let $F:\mu$-$\graph\ra\graph$ be the forgetful functor.  $F$-structures then do not extend along morphisms.  More plainly, if $H$ is a subgraph of a $\mu$-colorable graph $G$, then it is not necessarily the case that all $\mu$-colorings of $H$ extend to a $\mu$-coloring of $G$.

For a strongly compact cardinal $\mu$, it does hold that if all subgraphs of $G$ of lesser cardinality are $\mu$-colorable, then $G$ itself is $\mu$-colorable, regardless of the cardinality of $G$.  This example highlights that singular cardinal compactness is an expression of not so much the compactness as the approximability of a singular-sized structure by substructures \textsl{provided} those structures can be extended along ``underlying sets''.  The meaning of this is captured concisely by Def.~\ref{extend}.

In \cite{H}, Hodges gives an elegant axiomatization of singular compactness, based on set-systems that he calls ``algebras'' and ``bases''.  (See \cite{S} for Shelah's earlier unification of the cases he discovered.)  Let us show how Hodges's axioms are subsumed under the functorial formulation.  Roughly speaking, \cite{H} considers the situation $F:\A\ra\fB$ where $\A$ is the poset of free subalgebras of the target object $X$ and free inclusions among them, $\fB$ is the category of sets, and $F$ associates to a subalgebra its underlying set.  Hodges's setup follows, changing his notation and terminology slightly so as to make a better parallel with the rest of this paper.

Let $X$ be a set and $\mu$ an infinite cardinal.  Let $\cS$ be a set of subsets of $X$.  For each element $U$ of $\cS$, let a set $B(U)$, possibly empty, be given; each element of $B(U)$ is to be a subset of $\cS$.  ($X$ is to be thought of as the set underlying the target algebra; its size will be assumed singular later.  $\mu$ is a certain threshold cardinal.  The elements of $\cS$ are the sets underlying the subalgebras of $X$.  $B(U)$ is non-empty if and only if $U$ is a ``free algebra''.  In that case, each $\F\in B(U)$ should be thought of as data equivalent to giving a ``basis'' of $U$.  Intuitively, however, a basis $B$ of a free algebra is identified with the set of free algebras generated by subsets of $B$; in particular, a subset of $\cS$.)

The axioms are:

\begin{itemize}
\item[I.] $\cS$ is closed under unions of chains, and for every set $Y\subseteq X$ there exists $U\in\cS$ with $Y\subseteq U$ and $\card(U)\leqslant\card(Y)+\mu$.
\item[II.] For all $U\in\cS$ and $\F\in B(U)$: if $Z\in\F$ then $Z\subseteq U$; the collection $\F$ is closed under unions of chains; and for every set $Y\subseteq U$ there exists $Z\in\F$ with $Y\subseteq Z$ and $\card(Z)\leqslant\card(Y)+\mu$.
\item[III.] Given $\F\in B(U)$ and $V\in\F$, let us write $\F|_V$ for the subset $\{Z\in\F\;|\;Z\subseteq V\}$ of $\F$.  Now the axiom is: for all $\F\in B(U)$ and $V\in\F$, one has $\F|_V\in B(V)$.
\item[IV.] To ease notation, for $U, V\in\cS$, let us write $U\lhd V$ to mean that there exists $\F\in B(V)$ such that $U\in\F$.  (Intuitively, $U$ is a ``free factor'' of $V$.)  Now the axiom is: if $U\lhd V$ and $\F\in B(U)$, then there exists $\F'\in B(V)$ such that $\F'|_U=\F$.
\item[V.] Let $U_\alpha$, $\alpha\prec\kappa$, be a continuous chain of elements of $\cS$ and let $\F_\alpha\in B(U_\alpha)$ be such that $\F_\alpha = \F_\beta|_{U_\alpha}$ for all $\alpha\prec\beta\prec\kappa$.  Then
\[  \big\{ \cup_{\alpha\prec\kappa} Z_\alpha \;|\;
         Z_\alpha\in\F_\alpha\textup{\ where\ }Z_\alpha,\alpha\prec\kappa,\textup{\ is a continuous chain} \big\}
                                                                         \in B\big(\cup_{\alpha\prec\kappa} U_\alpha\big) \, . \]
\end{itemize}
For $U\in\cS$, let us write ``$U$ is free'' to mean ``$B(U)$ is non-empty''.  Singular compactness now takes the form: assume that a set-system satisfies Axioms I-V.  Assume $\mu<\card(X)$, that $\card(X)$ is singular, and that $U$ is free for all $U\in\cS$ with $\card(U)<\card(X)$.  Then $X$ is free.  (The text of \cite{H} omits the restriction `$\card(U)<\card(X)$,' but surely that is what's meant.)

To see this as a case of Theorem~\ref{main}, let $\A$ be the set of ``bases'', that is, subsets of $\cS$ of the form $\F$, with $\F\in B(U)$ for some $U\in\cS$.  Note that Axiom II implies that $\cup\F = U$ whenever $\F\in B(U)$; that is, a basis determines what it is a basis \textsl{of}.  Define the relation $\F\ra\F'$ to mean that $\F=\F'|_U$ where $U = \cup\,\F$.  It is immediate that $\ra$ is reflexive (and antisymmetric), and Axiom III implies that it is transitive, thus turning $\A$ into a small (a fortiori accessible)\ category.  Thanks to Axioms II and V, $\A$ has colimits of well-ordered chains; as is well-known, this implies that $\A$ has all directed colimits.  Let $F:\A\ra\Set$ send $\F$ to the set $\cup\,\F$.  Axioms II and V ensure that $F$ is a functor, preserving directed colimits.  The image of $F$ is a filter of $\Sub(X)$ that is $\mu$-dense by Axiom I and the assumption that all subalgebras of size less than $\card(X)$ are free.  Note that $U\subseteq V$ is the $F$-image of a morphism of $\A$ if and only if $U\lhd V$.  Axiom IV then says exactly that $F$-structures extend.  \qed

\section{Cellular version}
Some old examples of singular cardinal compactness due to Shelah, as well as recent ones involving filtered modules --- see~\cite{E} or~\cite{GT} --- while based on Hodges's axioms and thus part of the functorial formulation, are intricate enough to understand in detail.  In fact, the formulation of singular compactness in terms of cellular or `relatively free' objects sheds new light on many examples that we have already seen.

\begin{defn}  \label{cell-def}
Let $\X$ be a class of morphisms in a category $\fB$.  A morphism $m$ of $\fB$ is $\X$-\textsl{cellular} if there exist an ordinal $\alpha$ and a smooth chain $D:\big\{\beta\;|\;\beta\preceq\alpha\big\}\ra\fB$ such that $D(\beta)\ra D(\beta+1)$ is a pushout of a member of $\X$ whenever $\beta\prec\alpha$, and $m=D(\nil)\ra D(\alpha)$.  An object $X$ in $\fB$ is called $\X$-\textsl{cellular} if the unique morphism $0\to X$ from an initial object is $\X$-cellular.
\end{defn}

\begin{example}
Let $\fB$ be the category of topological spaces and let $\X=\{\partial D^n\ra D^n\;|\;n\in\N\}$, the inclusions of the boundary spheres of the unit balls in $\RR^n$.  $\X$-cellular maps are the \textsl{relatively cellular maps} and $\X$-cellular spaces are the \textsl{cellular spaces} of topology.
\end{example}

\begin{example}   \label{free}
Keeping the notation of the introduction, let $\X=\{F_\nil\ra F_\bullet\}$ in a category of algebras.  The $\X$-cellular objects are the free algebras.
\end{example}

\begin{example}
Let $U_i$, for $i\in I$, be a set of objects of $\fB$ and let $\X=\{0\ra U_i\;|\;i\in I\}$ where $0$ is an initial object.  The $\X$-cellular objects are the $\{U_i\}_{i\in I}$-decomposable ones, cf.\ Example~\ref{decomp}.
\end{example}

\begin{example}   \label{colourings}
Let $\fB$ be the category of graphs and graph homomorphisms.  Objects of this category are pairs $\lc U,A\rc$ where $U$ is a set and $A$ is a symmetric relation on $U$; morphisms are maps $u:U_1\to U_2$ preserving the relation.  This category is locally finitely presentable.  Let $\mu$ be a cardinal and $\X$ consist of all inclusions $\lc U,A\rc\to\lc U\cup\bullet,\;A\cup(U\times\bullet)\cup(\bullet\times U)\rc$ where $U$ has cardinality $<\mu$ and $\bullet$ is a singleton (we take a representative set of these morphisms).  Then $\X$-cellular objects are precisely graphs having coloring number $\leq\mu$, i.e.\ graphs whose vertices can be well-ordered such that property~(\ref{mu-short})\ holds for every vertex $v$.  Note that the well-ordering of vertices, which appeared somewhat artificial in Example~\ref{wellgraph}, is here a natural indicator of cellularity.
\end{example}

\begin{example}  \label{transversal}
Let $\fB$ be the category of (directed)\ bipartite graphs.  Objects of this category are triples $\lc U,V,E\rc$ where $U$ and $V$ are sets and $E\subseteq U\times V$.  Morphisms are $(u,v):\lc U_1,V_1,E_1\rc\to\lc U_2,V_2,E_2\rc$ where $u:U_1\to U_2$ and $v:V_1\to V_2$ are maps such that $(x,y)\in E_1$ implies that $(u(x),v(y))\in E_2$.  This category is locally finitely presentable.  Let $\X$ consist of the following three morphisms
\begin{enumerate}
\item[(1)] $f_1:\lc\nil,\nil,\nil\rc\to\lc\circ,\bullet,\star\rc$
\item[(2)] $f_2:\lc\nil,\nil,\nil\rc\to\lc\nil,\bullet,\nil\rc$
\item[(3)] $f_3:\lc\circ,\bullet,\nil\rc\to\lc\circ,\bullet,\star\rc$
\end{enumerate}
where $\nil$ is the empty set and $\circ$, $\bullet$, $\star$ are one-element sets.  Informally speaking, $f_1$ creates an edge, $f_2$ creates a vertex in the second partition, and $f_3$ creates an edge between two existing vertices.  $\X$-cellular objects are precisely bipartite graphs $\lc U,V,E\rc$ having a transversal, i.e.\ some injective map $t:U\ra V$ such that $(u,t(u))\in E$ for all $u\in U$.  That an $\X$-cellular graph possesses a transversal is immediate by transfinite induction.  For the converse, a graph $(U,V,E)$ with transversal $t$ can be built up from the empty graph with $U$ cells of type $f_1$; $V\setminus t(U)$ cells of type $f_2$; and $E\setminus\{(u,t(u))\;|\;u\in U\}$ cells of type $f_3$, in that order.
\end{example}

Much as free algebras extend along the functor $F:\Set_\mono\ra\Alg$ (but not, in general, along $F:\Set\ra\Alg$)\ cellular structures extend along suitable monomorphisms of posets indexing the cell attachments.  We recall some definitions from \cite{MRV}.

\begin{defn}
A poset $A$ is \textsl{good} if it is well-founded and has a least element.  A good poset $A$ is \textsl{$\mu$-good} if its initial segments $\{y\in A\;|\;y\preceq x\}$ have cardinality less than $\mu$ for all $x\in A$.  An element $x\in A$ is \textsl{limit} if the strict initial segment $\{y\in A\;|\;y\prec x\}$ is non-empty and does not have a top element.  A diagram $D:A\ra\fB$ is \textsl{smooth} if for every limit $x\in A$, the diagram $D(-,x): Dy\ra Dx|_{y\prec x}$ is a colimit cocone on the restriction of $D$ to $\{y\in A\;|\;y\prec x\}$.
\end{defn}

\begin{prop}   \label{cell-extend}
Let $\fB$ be a cocomplete category, $\mu$ a regular uncountable cardinal and $\X$ a class of morphisms in $\fB$ with $\mu$-presentable domains.  Let $\A$ be the category whose objects are smooth diagrams $D:A\to\fB$ where
\begin{itemize}
\item[-] $A$ is a $\mu$-good poset
\item[-] for all $a\in A$, the object $D(a)$ is $\X$-cellular and $\mu$-presentable in $\fB$
\item[-] for all $a\preceq b\in A$, the morphism $D(a)\to D(b)$ is $\X$-cellular.
\end{itemize}
A morphism from $D_1:A_1\to\fB$ to $D_2:A_2\to\fB$ is a monotone map $m:A_1\ra A_2$ such that $D_1=D_2m$ and $m$ is a downward-closed embedding, i.e.\ $m(a)\preceq m(b)$ if and only if $a\preceq b$ and if $b\preceq m(a)$ for $a\in D_1$, $b\in D_2$ then $b$ is in the image of $m$.

Let $F:\A\to\fB$ be the functor sending the diagram $D:A\ra\fB$ to its colimit.  Then $F$-structures extend along morphisms.
\end{prop}

\begin{proof}
We build heavily on the results of \cite{MRV}.  By \cite{MRV}~4.14 and 4.15(1), an object of $\fB$ is $\X$-cellular if and only if it belongs to the image of $F$.  Let $D_1:A_1\to\fB$ and $D_2:A_2\to\fB$ be objects of $\A$ and $m$ a morphism of $\A$ determined by a downward-closed embedding $m:A_1\to A_2$.  Let $F(m):F(D_1)\to F(D_2)$ be the induced morphism on colimits.  Let $E:B\to\fB$ be another object of $\A$ with $F(E)\cong F(D_1)$.  By \cite{MRV}~4.10, there is an extension $E^\ast:B^\ast\to\fB$ of $E$ to a $\mu$-good $\mu$-directed diagram of $\X$-cellular morphisms such that $B$ is downward-closed in $B^\ast$ and the induced morphism $F(E)\to F(E^\ast)$ on colimits is an isomorphism.  It immediately follows from the construction in \cite{MRV}~4.10 that $E^\ast$ is a diagram of $\X$-cellular $\mu$-presentable objects, that is, an object of $\A$.  Now, we have an $\X$-cellular morphism $F(m):F(D_1)\to F(D_2)$ and an isomorphism $F(D_1)\to F(E^\ast)$ and we need to make $E^\ast$ to be an initial segment in a $\mu$-good diagram $D^\prime:A^\prime\to\fB$ of $\X$-cellular $\mu$-presentable objects and $\X$-cellular morphisms such that $F(D^\prime)\cong F(D_2)$ and the composition $E\to E^\ast\to D^\prime$ yields the desired reparametrization of $m$.  Here we can imitate the proof of \cite{MRV}~4.11.  Let $m_{01}$ be the first morphism in the cellular decomposition of $m$.  Then we have a pushout
\[  \xymatrix@=3pc{
A_1 \ar[r]^{m_{01}} & A_{11} \\
X \ar [u]^{u_1} \ar [r]_{g} &
Y \ar[u]_{u_2}  }  \]
where $g\in\X$.  Thus the composition of $u_1$ with the isomorphism $F(D_1)\to F(E^\ast)$ is a morphism from $X$ to the colimit of $E^\ast$.  Since $E^\ast$ is $\mu$-directed and $X$ is $\mu$-presentable, this composition factorizes through a morphism $u:X\to E^\ast(b)$ for some $b\in B^\ast$.  We take all pushouts of $m_{01}$ along $E^\ast(t)u$ where $t:b\to b'$ is a morphism in $B^\ast$. We get a new $\mu$-good $\mu$-directed diagram $E_1^\ast:B_1^\ast\to\fB$ of $\X$-cellular $\mu$-presentable objects and a downward-closed embedding $B^\ast\to B_1^\ast$.  Moreover $F(E_1^\ast)\cong A_{11}$.  Then $D'$ is obtained by continuing this procedure where, in limit
steps, we take the $\ast$-closure of the union of preceding diagrams.
\end{proof}

\begin{thm}   \label{cellular}
Let $\fB$ be a locally finitely presentable category, $\mu$ a regular uncountable cardinal and $\X$ a class of morphisms with $\mu$-presentable domains.  Let $X\in\fB$ be an object with $\max\{\mu,\card(\fp\fB)\} <\|X\|$.  Assume
\begin{itemize}
\item[(i)] $\|X\|$ is a singular cardinal
\item[(ii)] there exists $\phi < \|X\|$ such that for all regular $\kappa$ with $\phi < \kappa < \|X\|$, there exists a dense $\kappa$-filter of $\Sub(X)$ consisting of $\X$-cellular objects.
\end{itemize}
Then $X$ is $\X$-cellular.
\end{thm}

\begin{proof}
We retain the notation of Prop.~\ref{cell-extend}.  The category $\A$ has directed colimits, computed as in the category of posets (i.e.\ on underlying sets).  Since initial segments of $\mu$-good posets have cardinality less than $\mu$, any object $D:A\ra\fB$ of $\A$ is the $\mu$-directed union of the restrictions of $D$ to downward closed subposets of $A$ of cardinality less than $\mu$.  Since objects $A\ra\fB$ of $\A$ with $\card(A)<\mu$ are $\mu$-presentable and there is only a set of isomorphism classes of such objects, $\A$ is a $\mu$-accessible category with filtered colimits.

We are going to apply Theorem~\ref{main} to the functor $F:\A\to\fB$ sending the diagram $D:A\ra\fB$ to its colimit.  $F$ obviously preserves filtered colimits.  Prop.~\ref{unif} implies that $\mu_F\leqslant\mu$ and Prop.~\ref{cell-extend} verifies the last condition of \ref{main}.
\end{proof}

Examples \ref{free} through \ref{colourings} recover familiar cases of singular compactness.  Example~\ref{transversal} yields Shelah's theorem on set transversals.  To any set of sets $\cS$, one can associate the bipartite graph $\lc\cS,\cup\,\cS,\in\rc$ where there is an edge from $A\in\cS$ to $a\in\cup\,\cS$ if and only if $a\in A$.  A transversal of $\cS$, i.e.\ a set of disjoint representatives, is then the same as a transversal of the associated bipartite graph.  Note that a subobject $\lc U_0,V_0,E_0 \rc$ of $\lc U,V,E \rc$ is given by injective maps $i: U_0\ra U$ and $j: V_0\ra V$ that preserve edges.  In particular, every subset $\cS_0\subset\cS$ gives rise to a sub-bipartite graph, but not all sub-bipartite graphs of $\cS$ arise this way.

\begin{prop}
Let $\lambda_0 < \lambda$ be infinite cardinals with $\lambda$ singular.  Let $\cS$ be a set of sets such that $\card(\cS)=\lambda$ and $\card(A)\leq\lambda_0$ for all $A\in\cS$.  Suppose that $\cS_0$ has a transversal for all subsets $\cS_0\subset\cS$ with $\card(\cS_0) < \lambda$.  Then $\cS$ itself has a transversal.
\end{prop}

\begin{proof}
Apply Theorem~\ref{cellular} to the object $X=\cS$ with the category $\fB$ played by bipartite graphs, the morphisms $\X$ of Example~\ref{transversal}, and $\mu=\aleph_1$.  Note that for a bipartite graph $\lc U,V,E\rc$,
\[   \|\,\lc U,V,E\rc\,\| = \max\{\card(U),\card(V)\}  \]
so $\|\cS_0\| = \card(\cS_0)$ as long as $\card(\cS_0)\geq\lambda_0$.  (We only consider infinite graphs, and do not notationally distinguish sets-of-sets from graphs.)  In particular, $\|\cS\|=\lambda$.  Consider the set $\F$ of all subsets $\cS_0\subset\cS$ with $\card(\cS_0) < \lambda$.  We claim $\F$ is a $\kappa$-filter of $\Sub(\cS)$ for every regular $\kappa < \lambda$.  Indeed, the colimit (i.e.\ union)\ of a $\leq\kappa$-long chain of $\kappa$-presentable objects in $\F$, belongs to $\F$.  Moreover, $\F$ is $\kappa$-dense in $\Sub(\cS)$ for any $\lambda_0 < \kappa < \lambda$.  Indeed, any subobject $\lc U,V,E \rc$ of $\cS$ is contained in $\lc U,\cup\,U,\in\rc$ and for any $\kappa > \lambda_0$, if $\lc U,V,E \rc$ is $\kappa$-presentable, so is $\lc U,\cup\,U,\in\rc$.  This verifies condition (ii)\ of \ref{cellular}.  (Note that it would suffice to assume that there exists $\phi < \lambda$ so that all subsets $\cS_0\subset\cS$ with $\phi < \card(\cS_0) < \lambda$ have a transversal.)
\end{proof}

Theorem~\ref{cellular} permits the collection of morphisms $\X$ to be a proper class, as long as the ranks of their domains are bounded by some cardinal $\kappa$.  For a given object $X$, only a subset $\X_0$ of the class $\X$ is needed to exhibit a dense filter of subobjects of $X$ as $\X$-cellular, i.e.\ to satisfy condition (ii)\ of Theorem~\ref{cellular}.  Obviously, $X$ will already be $\X_0$-cellular, but the set $\X_0$ may depend on $X$.  There are also cases when a subset $\X_0\subset\X$ can be shown to exist such that the classes of $\X$-cellular and $\X_0$-cellular objects are the same.  This happens in deep work of Eklof et al.\ that we alluded to above.

\begin{example}\label{filtered}
Let $\fB$ be the category of $R$-modules and $\cQ$ a set of $R$-modules.  Let $\X$ be the class of all $\cQ$-monomorphisms, i.e.,
monomorphisms whose cokernel is isomorphic to a member of $\cQ$.  Then $\X$-cellular objects are precisely the $\cQ$-filtered modules.  By a result of Saor\'{\i}n and \v{S}\v{t}ov\'{\i}\v{c}ek~\cite{SS}, $\X$-cellular objects coincide with $\X_0$-cellular ones for a certain subset $\X_0$ of $\X$.  Thus Theorem~\ref{cellular} implies the singular compactness theorem for $\cQ$-filtered modules.  (See \cite{E} or \cite{GT}.)  The threshold cardinal $\max\{\kappa,\card(\fp\fB)\}$ of \ref{cellular} works out to be an uncountable upper bound for the presentability ranks of modules from $\cQ$.

If we take $\fB$ to be the category of abelian groups and $\cQ$ to consist of the group of integers, $\X$-cellular objects are precisely
free abelian groups and we get Shelah's~\cite{S} original singular compactness for free abelian groups again.
\end{example}

\begin{remark}
In many examples of the cellular version of singular compactness, the generating cellular maps in $\X$ are monomorphisms, and monomorphisms are preserved by pushouts in $\fB$.  In particular, cellular subobjects of the target $X$ of Theorem~\ref{cellular} will be built up as smooth chains of subobjects of $X$.  But it is easy to find locally finitely presentable categories where the pushout of a monomorphism need not be a monomorphism.  For example, in the category of rings and homomorphisms, the below square is a pushout:
\[  \xymatrix{  \ZZ \ar[r]\ar[d] & R \ar[d] \\
                \QQ \ar[r] & \te  }  \]
where $R$ is any finite ring, and $\te$ is the terminal (one-element)\ ring.  A subobject of $X$ thus may be cellular, with this fact not being ``witnessed'' by any smooth chain that lies entirely within the subobject lattice of $X$.  Perhaps this is the most significant departure of Theorem~\ref{cellular} from singular cardinal compactness as applying to set-systems and inclusions.
\end{remark}

\end{document}